\documentclass[12pt]{amsart}
\usepackage{amscd,amssymb,graphics,color,a4wide,hyperref,verbatim}
\usepackage{graphicx}
\usepackage{psfrag} 
\usepackage{stmaryrd} 

\usepackage{mathrsfs}
\input xy
\xyoption{all}

\newcommand{\op}[1]{\operatorname{#1}}

\newcommand{\Cone}{\operatorname{{Cone}}}

\DeclareFontFamily{U}{rsf}{}
\DeclareFontShape{U}{rsf}{m}{n}{
  <5> <6> rsfs5 <7> <8> <9> rsfs7 <10->  rsfs10}{}
\DeclareMathAlphabet{\mathscr}{U}{rsf}{m}{n}

\newtheorem{theorem}{Theorem}[section]

\newtheorem{proposition}[theorem]{Proposition}

\theoremstyle{definition}
\newtheorem{definition}[theorem]{Definition}

\newtheorem{example}[theorem]{Example}

\newtheorem*{acknowledgement}{Acknowledgement}

\theoremstyle{remark}

\numberwithin{equation}{section}

\newcommand{\NN} {\mathbb{N}}
\newcommand{\ZZ} {\mathbb{Z}}
\newcommand{\QQ} {\mathbb{Q}}
\newcommand{\RR} {\mathbb{R}}
\newcommand{\CC} {\mathbb{C}}

\newcommand{\HH} {\mathbb{H}}
\newcommand{\PP} {\mathbb{P}}
\renewcommand{\AA} {\mathbb{A}}

\newcommand {\shF} {\mathcal{F}}

\newcommand {\shL} {\mathcal{L}}

\newcommand {\shO} {\mathcal{O}}

\newcommand {\shT} {\mathcal{T}}

\newcommand {\shX} {\mathcal{X}}

\newcommand {\Bl} {\operatorname{Bl}}

\newcommand {\C} {\mathbb{C}}

\newcommand {\dlog} {\operatorname{dlog}}

\newcommand {\dual} {{\vee}}

\newcommand {\End} {\operatorname{End}}

\newcommand {\GL} {\operatorname{GL}}

\newcommand {\Gr} {\operatorname{Gr}}

\newcommand {\Hom} {\operatorname{Hom}}

\newcommand {\id} {\operatorname{id}}
\newcommand {\im} {\operatorname{im}}
\renewcommand{\Im} {\operatorname{Im}}

\newcommand {\Int} {\operatorname{Int}}

\newcommand {\lra} {\longrightarrow}

\newcommand {\new} {{\operatorname{new}}}

\renewcommand{\P} {\mathscr{P}}

\newcommand {\Proj} {\operatorname{Proj}}

\newcommand {\ra} {\to}

\newcommand {\Sing} {\operatorname{Sing}}

\newcommand {\Spec} {\operatorname{Spec}}

\newcommand {\sra}{\twoheadrightarrow}

\newcommand {\Sym} {\operatorname{Sym}}

\newcommand {\Tot} {\operatorname{Tot}}

\def\mydate{\ifcase\month \or January\or February\or March\or
April\or May\or June\or July\or August\or September\or October\or 
November\or December\fi \space\number\day,\space\number\year}

\newlength{\picwidth} 
\newlength{\miniwidth} 
\newlength{\nanowidth} 
\newlength{\melowidth} 
\newlength{\leftminiwidth} 
\newlength{\rightminiwidth} 
\newlength{\minipagewidth} 

\begin{document}
\def\mapright#1{\smash{
 \mathop{\longrightarrow}\limits^{#1}}}
\def\mapleft#1{\smash{
 \mathop{\longleftarrow}\limits^{#1}}}
\def\exact#1#2#3{0\to#1\to#2\to#3\to0}
\def\mapup#1{\Big\uparrow
  \rlap{$\vcenter{\hbox{$\scriptstyle#1$}}$}}
\def\mapdown#1{\Big\downarrow
  \rlap{$\vcenter{\hbox{$\scriptstyle#1$}}$}}
\def\dual#1{{#1}^{\scriptscriptstyle \vee}}
\def\invlim{\mathop{\rm lim}\limits_{\longleftarrow}}
\def\rto{\raise.5ex\hbox{$\scriptscriptstyle ---\!\!\!>$}}

\input epsf.tex
\title
[Mirror duality of Landau-Ginzburg models via Discrete Legendre Transforms]
{Mirror duality of Landau-Ginzburg models via Discrete Legendre Transforms}
\author{Helge Ruddat}

\address{JGU Mainz, Institut f\"ur Mathematik, Staudingerweg 9, 55099 Mainz}
\email{ruddat@uni-mainz.de}

\begin{abstract}
We recall the semi-flat Strominger-Yau-Zaslow (SYZ) picture of mirror symmetry and discuss the transition from the Legendre transform to a discrete Legendre transform in the large complex structure limit. We recall the reconstruction problem of the singular Calabi-Yau fibres associated to a tropical manifold and review its solution in the toric setting. We discuss the monomial-divisor correspondence for discrete Legendre duals and use this to give a mirror duality for Landau Ginzburg models motivated from the SYZ perspective and Floer theory. We mention its application for the construction of mirror symmetry partners for varieties of general type and discuss the straightening of the boundary of a tropical manifold corresponding to a smoothing of the divisor in the complement of a special Lagrangian fibration.
\end{abstract}

\maketitle
\setcounter{tocdepth}{1}
\tableofcontents
\bigskip

\section{Strominger-Yau-Zaslow fibrations and the mirror of $(\CC^*)^n$}
\label{section1}
We give a summary of the semi-flat picture of mirror symmetry following \cite[\S6-8]{Clay09} and discuss the example of an algebraic torus. Further references for the material are 
\cite{Mi04}, \cite{Le05}, \cite{CM06}, \cite{Au07}, \cite{Gr08}, \cite{CL08}, \cite{CLL10} and most recently \cite{Gr12}.
Hitchin \cite{Hi97} first noticed the importance of the Legendre transform in this context.
A Legendre transform already appeared in \cite{Gu94} in a closely related context without the awareness of mirror symmetry and special Lagrangians.

Mirror symmetry has become intrinsic to the Calabi-Yau geometry by the work of Strominger-Yau-Zaslow \cite{SYZ96} (short: SYZ), suggesting to explain the mirror duality of two Calabi-Yau manifolds $X$, $\check X$ as a duality of torus fibrations. 
There are supposed to be $C^\infty$-maps 
$$f:X\ra B,\quad  \check f:\check X\ra B$$ 
with fibres homeomorphic to $(S^1)^n$ for $n=\dim_\CC X=\dim_\RR B$, in fact if $f^{-1}(b)=V/\Lambda$ for a real vector space $V$ with lattice $\Lambda\cong\ZZ^n$ then $\check f^{-1}(b)=V^*/\Lambda^*$ where $V^*=\Hom(V,\RR), \Lambda^*=\Hom(\Lambda,\ZZ)$. 
Moreover, in the strong form of Strominger-Yau-Zaslow, the fibres of $f$ and $\check f$ are required to be special Lagrangian, so by definition the restriction to the fibres of the symplectic form $\omega$ and the imaginary part of a fixed holomorphic volume form $\Omega$ vanish respectively.
The base $B$ carries the structure of a real affine manifold in two ways as follows. 
The transitions between coordinate charts of $B$ are going to be elements of $\GL_n(\ZZ)\ltimes\RR^n$ respectively. 

One affine structure is determined by the complex structure of $X$ and alternatively also by the symplectic structure of $\check X$.
The other affine structure is determined by the symplectic structure of $X$ and alternatively also by the complex structure of $\check X$.
Let $\nu$ denote the vector field on $f^{-1}(b)$ given as a lift of a tangent vector $\bar\nu$ at a point $b\in B$ then the contraction of $\omega$ (respectively $\im\Omega$) by $\nu$ yields a one-form (respectively $(n-1)$-form) on $f^{-1}(b)$. That these are independent of the lift chosen follows from $f^{-1}(b)$ being special Lagrangian.
McLean showed (\cite[\S6.1]{Clay09}) that these two forms on $f^{-1}(b)$ are both closed if and only if the infinitesimal deformation $\bar\nu$ of $f^{-1}(b)$ preserves the special Lagrangian property (which is true for a special Lagrangian fibration). Moreover, these two forms can be shown to be Hodge-star dual on $f^{-1}(b)$, so first order Lagrangian deformations correspond to harmonic one-forms on the Lagrangian.
McLean proves that the moduli space of special Lagrangians is unobstructed \cite[Thm 3--4]{ML98}.
One deduces from this that $B$ is locally the moduli space of the fibres of $f$ as well as $\check f$.
The just constructed maps descend to isomorphisms on cohomology
\begin{equation}
\label{flatstructures}
\begin{array}{c}
\shT_{B,b}\cong_\omega H^1(f^{-1}(b),\RR),\\[2mm]
\shT_{B,b}\cong_{\im\Omega} H^{n-1}(f^{-1}(b),\RR),
\end{array}
\end{equation}
which give the tangent bundle two usually different flat connections. 
To distinguish the two, we denote the manifold $B$ with the flat structure coming from $f$ and $\omega$ by $\check B$ whereas the manifold with flat structure derived from $f$ and $\im\Omega$ keeps the name $B$. For either of these, we call a set of coordinates $\{y_j\}$ \emph{affine} if $\partial_{y_j}$ are flat with respect to the respective flat structure.
We also obtain a local systems of integral tangent vectors $\Lambda_B\subset \shT_B$ isomorphic to the integral cohomology $H^{n-1}(f^{-1}(b),\ZZ)\subset H^{n-1}(f^{-1}(b),\RR)$ and 
similarly a system $\Lambda_{\check B}\subset \shT_{\check B}$. 
A set of coordinates $\partial_{y_j}$ on $B$ (resp. $\check B$) is called \emph{integral affine} if $\partial_{y_j}\in\Lambda_B$ (resp. in $\partial_{y_j}\in\Lambda_{\check B}$) and they form a basis over $\ZZ$. Thus, $B$ and $\check B$ are real affine manifolds with coordinate transitions in $\GL_n(\ZZ)\ltimes\RR^n$.

We assume that the torus bundle $f$ is oriented and obtain from the second equation in \eqref{flatstructures}, $\shT_{B,b}\cong (H^1(f^{-1}(b),\RR))^*=H_1(f^{-1}(b),\RR)$. 
Under this isomorphism, $\Lambda_B$ becomes $H_1(f^{-1}(b),\ZZ)$, so we have $X\cong \shT_B/\Lambda_B$ as topological manifolds.
Alternatively, we may also use $\shT^*_{\check B}\cong H_1(f^{-1}(b),\RR)$ by means of the first equation in \eqref{flatstructures} to reconstruct $X$. We summarize
\begin{equation}
\label{topreconstruct}  
\shT^*_{\check B}/\Lambda^*_{\check B} \ \underset{\omega}{\cong} \ X\ \underset{\im\Omega}{\cong}\ \shT_B/\Lambda_B.
\end{equation}

\begin{figure}
\resizebox{0.6\textwidth}{!}{
\input{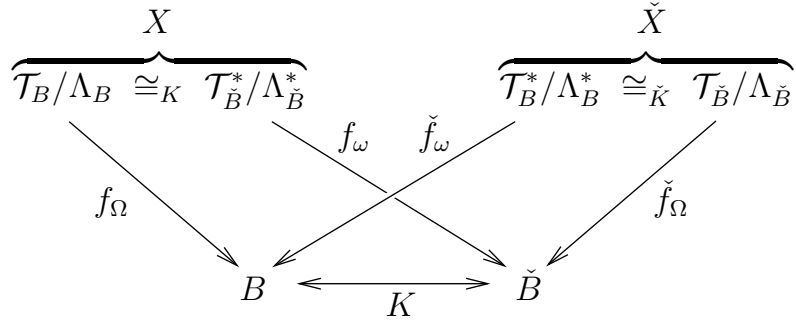}
}
\caption{$f_\Omega,\check f_\Omega$ are logarithm maps, $f_\omega, \check f_\omega$ are moment maps.}
\label{logandmoment}
\end{figure}

We can play the same game with $\check f:\check X\ra B$ in place of $f:X\ra B$ and the definition of $SYZ$ mirror duality for $X,\check X$ is the statement that this is supposed to yield identical affine manifolds $B$, $\check B$ with swapped roles, i.e. the flat structure on $B$ derives from the symplectic structure on $\check X$ and the flat structure on $\check B$ from the holomorphic structure on $\check X$, see Fig.~\ref{logandmoment}.

The work of Gross and Siebert on mirror symmetry by means of toric degenerations, starting out with \cite{GS03}, was motivated by reverse engineering $X$ and $\check X$ from $B$.
The real difficulty arises when $X,\check X$ are intended to be compact since then $f,\check f$ need to have singular fibres and the affine structures need to have singularities as well.
We will not deal with singularities before \S\ref{MirdualLGmodels} but we adopt the point of view of reconstructing $X$ and $\check X$ from $B$.
In what we discussed so far, at least topologically by \eqref{topreconstruct}, the reconstruction of $X,\check X$ is straightforward once we know $\Lambda_B$ and $\Lambda_{\check B}$. In fact, this is a datum we need to fix in addition to $B$ and $\check B$.
This topological picture can be enhanced as follows. 
Given the real affine manifold $B$, we have 
\begin{itemize}
\item[(A)]\quad a canonical symplectic structure on $\check X:=\shT^*_B/\Lambda^*_B$ locally given by
$\omega=\sum_j d\bar{x}_j\wedge dy_j$ where $y_j$ are affine coordinates of $B$ and $\bar{x}_j=\partial_{y_j}$,
\item[(B)]\quad a canonical complex structure on $X:=\shT_B/\Lambda_B$ locally given by
complex coordinates $z_j=x_j+iy_j$ where $y_j$ are integral affine coordinates of $B$, $x_j=dy_j$ and $i=\sqrt{-1}$. The holomorphic volume form is $\Omega=dz_1\wedge\dots\wedge dz_n$.
We set $w_j=e^{2\pi i z_j}$.
\end{itemize}
Note that integrality of the coordinates only matters in (B).
To obtain the complementary parts, i.e., the symplectic structure on $X$ and complex structure on $\check X$, one uses the structure of a Hessian metric $g$ on $B$. 
We obtain K\"ahler manifold $X$ and $\check X$ by applying (A), (B) on the respective dual side using $g$ to identify the tangent and cotangent bundle of $B$.
More explicitly, $g$ is locally given as $g_{ij}=\partial_{y_i}\partial_{y_j}K$ for some smooth strictly convex function $K:B\ra\RR$. 
Mirror duality appears in this setup in the disguise of the Legendre transform, see \cite[Prop. 6.4]{Clay09}:
\begin{definition} 
Given a real affine manifold $B$ with Hessian metric $g$, the \emph{Legendre transform} is the real affine manifold $\check B$ which is homeomorphic to $B$ with coordinates given by
$\check y_j:=\partial_{y_j} K$ 
(where $y_j$ are local affine coordinates on $B$ and $K$ is a local potential defining $g$) and dual potential $\check K:\check B\ra\RR$,
$$\check K(\check y_1,\dots,\check y_n)=\sum_j \check y_j y_j - K(y_1,\dots,y_n).$$
\end{definition}
Note that also the integral structure dualizes: dual integral affine coordinates are those that are the Legendre dual of integral affine coordinates.

The symplectic structure on $X$ and the complex structure on $\check X$ is given directly by 
\begin{equation} \label{dualviaK}
\begin{array}{rcl}
\omega&=&2i\partial\bar\partial (K\circ f)=\frac{i}2\sum g_{jk}dz_j\wedge d\bar z_k,\\
\bar z_j&=&\bar x_j+i\partial_{y_j} K,
\end{array}
\end{equation}
see \cite[Prop 3.2]{Gr08}, \cite[Prop. 6.15]{Clay09}.

The manifold $X$ (resp $\check X$) is Ricci-flat (i.e., $\omega^n=c\Omega\wedge\bar\Omega$ for some $c\in\CC$) if and only if $\det(\partial_{y_i}\partial_{y_j}K)=\det(g)$ is constant as follows from \eqref{dualviaK}.

We dicuss the following integrated version of the two affine structures which was pointed out to the author by Denis Auroux. 
It gives a hint at why mirror symmetry would exchanges periods and Gromov-Witten-invariants. Moreover, it leads towards Landau-Ginzburg potentials.
Let $X$ be a Calabi-Yau with K\"ahler form $\omega$ and non-vanishing holomorphic volume form $\Omega$. The affine manifold $B$ is the moduli space of special Lagrangian tori in $X$, i.e., the moduli of manifolds $L$ homeomorphic to $(S^1)^n$ with $\omega|_L=0$ and $\Im\Omega|_L=0$ (more generally one allows for a phase $\theta\in\RR$, i.e., $\Im(e^{i\theta}\Omega)|_L=0$). Moreover, $\check X$ is given as the moduli space of pairs $(L,\nabla)$ where $L$ is special Lagrangian and $\nabla$ is a flat $U(1)$-connection of the trivial bundle with fibre $\CC$ on $L$. 
The information of $\nabla$ is equivalent to a map of groups $H_1(L,\ZZ)\ra U(1)$. \\[-2mm]

\begin{minipage}[c]{0.42\textwidth}
\begin{center}
\input{tori.pstex_t}
\end{center}
\end{minipage}\qquad\qquad\qquad
\begin{minipage}[c]{0.4\textwidth} 
The local integral affine coordinates on the base are then given as
\begin{equation} 
\label{integrals}
\begin{array}{rcl}
y_i&=&\int_{\Gamma_i}\omega, \\ 
\check y_i&=&\int_{\Gamma^*_i}\Im\Omega
\end{array}
\end{equation}
where $\Gamma_i\in H_2(X,L\cup L')$ are cylinders traced out by a basis $\{\gamma_i\}$ of $H_1(L,\ZZ)$ as we move $L$ to $L'$ and $\Gamma^*_i\in H_n(X,L\cup L')$ are traced out by a basis $\{\gamma^*_i\}$ of $H_{n-1}(L,\ZZ)$ as we move $L$ to $L'$.
\end{minipage}

\begin{example}[The mirror dual of $(\CC^*)^n$]
\label{Cstarn}
The simplest example is 
$X=(\CC^*)^n$. Its complex structure is indeed given as in (B) if we identify
$B=\RR^n$, $\shT_B=\RR^n\times\RR^n$, $\Lambda=\ZZ^n$ where the latter is naturally contained in the second factor of $\shT_B$. On the universal covers of $(\CC^*)^n$ and $\shT_B/\Lambda$ we set $z_j=x_j+iy_j$ where $z_j$ are standard coordinates on $\CC^n$, $y_j$ are standard coordinates on $B$, $x_j=dy_j$ and
$w_j=e^{2\pi i z_j}$ are standard coordinates on $(\CC^*)^n$. We thus obtain
$$f:(\CC^*)^n\ra B,\qquad (w_1,\dots,w_n)\mapsto \frac{-1}{2\pi}(\log|w_1|,\dots,\log|w_n|)=(y_1,...,y_n).$$
The holomorphic volume from is given by (B) as follows, we additionally pick the following symplectic form 
$$
\begin{array}{rcccl}
\Omega&=&\frac1{(2\pi i)^{n}}\dlog w_1\wedge\dots\wedge \dlog w_n &=& dz_1\wedge\dots\wedge dz_n,\\  
\omega&=&\frac{-1}{(2\pi)^{2}}\sum_j \dlog r_j \wedge d\theta_j&=&\sum_j dx_j\wedge dy_j
\end{array}
$$
where $w_j=r_j e^{i\theta_j}$. This choice turns $f$ into a special Lagrangian fibration with $y_j=\check y_j$ as follows directy from \eqref{integrals}. It determines $K=\frac12\sum y_j^2$ up to a constant and $g$ is the standard metric on $B$.
We conclude from $y_j=\check y_j$ and $\shT_B\cong \shT^*_B$, that
\begin{center}
\fbox{
\emph{the SYZ mirror dual of $((\CC^*)^n,\Omega,\omega)$ is $((\CC^*)^n,\Omega,\omega)$.}}
\end{center}
\end{example}

The setup in this example is very special in the sense that the two sets of affine coordinates on $B$ coincide. It is easy to check that indeed $\omega=2i\partial\bar\partial(K\circ f)$.
More generally, the situation can be diagrammed as in Figure~\ref{logandmoment}.
As verified in \cite{Au07},\,Prop.\,4.2, $\check f_\omega$ coincides with the moment map associated to $\omega$ and the natural fibrewise $(S^1)^n$-action on $\check X=\shT^*_B/\Lambda^*$ and similarly for $f_\omega$. Moreover, as in the above example, $f_\Omega$ is expressible as the map $\frac{-1}{2\pi}\log|\cdot|$ componentwise in the complex coordinates $w_j$ on $X=\shT_B/\Lambda$.

\begin{figure} 
\includegraphics[width=0.48\textwidth]{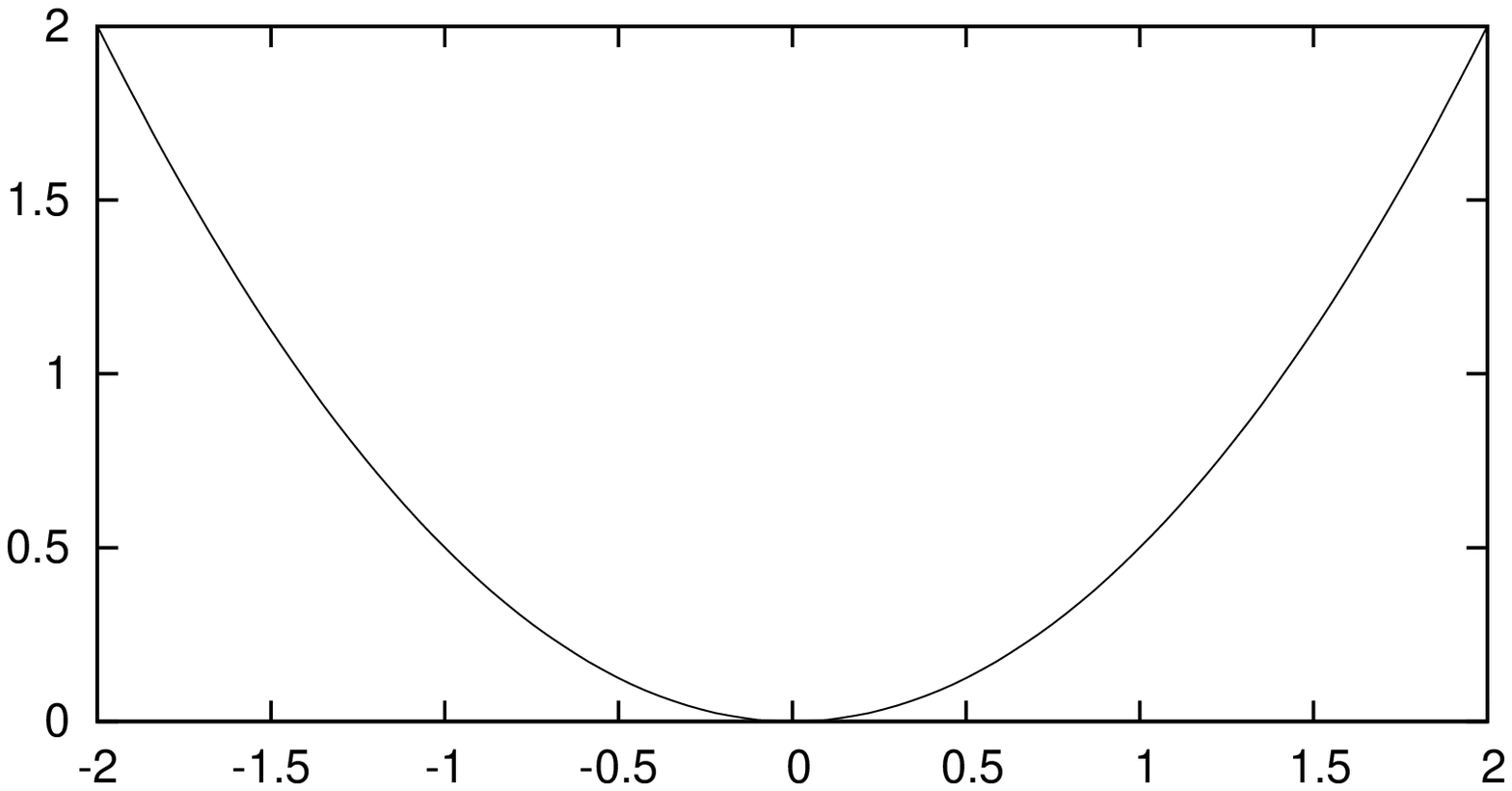}
\includegraphics[width=0.48\textwidth]{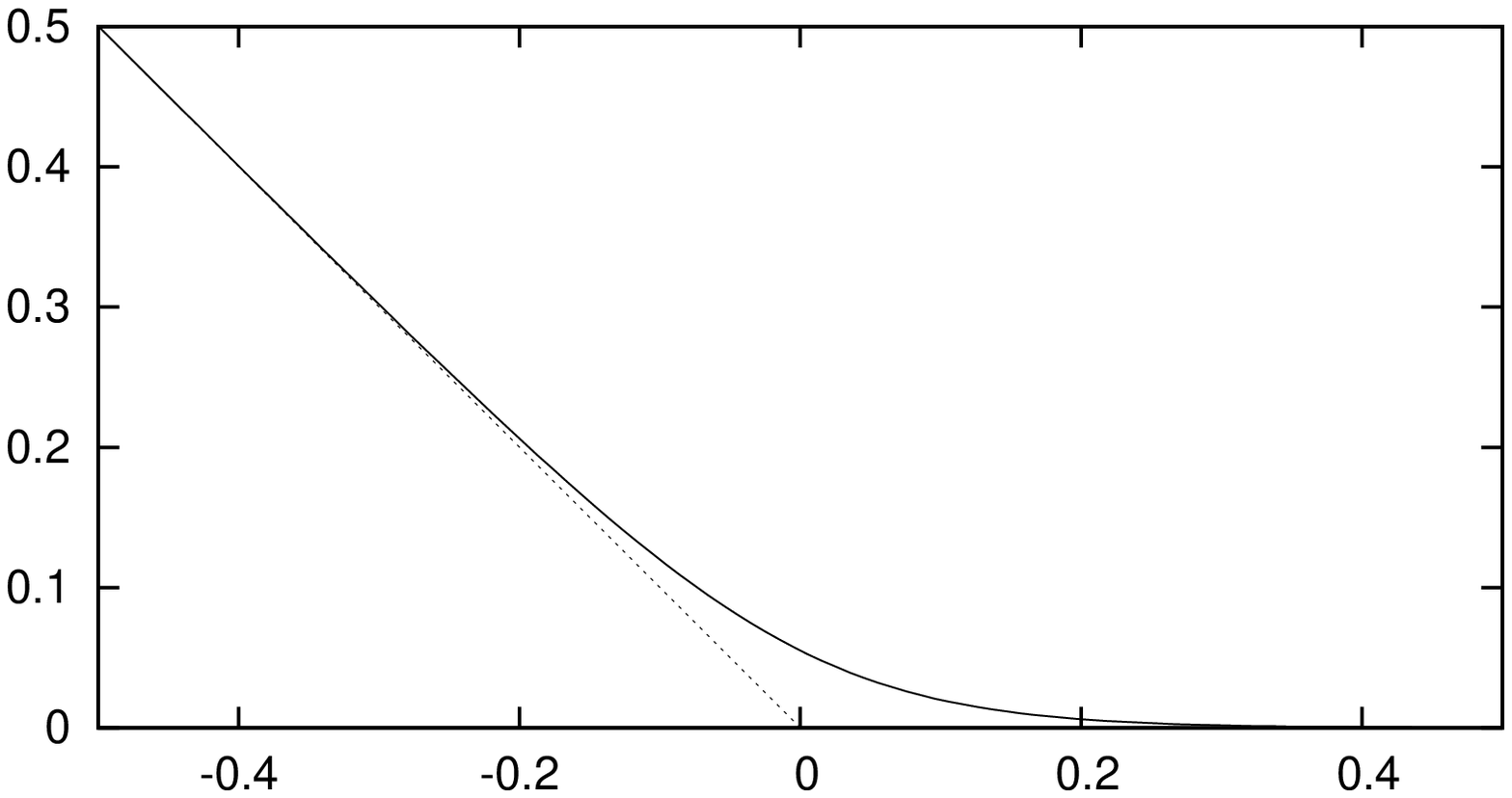}
\caption{$K(y)=\frac12y^2$ and $K(y)=\frac1{4\pi} \op{log}(1 + e^{-4\pi y})=-\int_{-\infty}^y\frac{e^{-4\pi u}}{1+e^{-4\pi u}}du$}
\label{figpotentials}
\end{figure}

\begin{example}[Further mirror duals of $(\CC^*)^n$]
\label{furtherCstarn}
While there aren't any interesting alternative algebraic choices for $\Omega$ in the previous example, there is a variety of choices for $\omega$: for each equivariant embedding 
$$\varphi:(\CC^*)^n\ra (\CC^*)^{m+1}/\CC^*,\qquad (w_1,\dots,w_n)\mapsto 
(\textstyle\prod_{k=1}^n w_k^{a_{0k}}:\dots:\textstyle\prod_{k=1}^n w_k^{a_{mk}})$$
we can take $\omega=\varphi^*\omega_{\op{FS}}$ where $\omega_{\op{FS}}$ is the Fubini-Study form on $\PP^m=(\CC^{m+1}\setminus\{0\})/\CC^*$ (normalized by $\int_{\PP^1}\omega=1$), 
i.e., 
\begin{equation} \label{FSform}
\pi^*\omega_{\op{FS}}=\frac{i}{2\pi}\partial\bar\partial\log\|z\|^2
\end{equation}
for $\pi:\CC^{m+1}\setminus\{0\}\ra\PP^m$ the natural projection. We want to compute $f_\omega$.
Let $S^{2m+1}=\{z\,|\,\|z\|=1\}$ denote the unit sphere in $\CC^{m+1}$. A straightforward computation shows that
\begin{equation} \label{FSisstandard}
\left.\left(\partial\bar\partial\log\|z\|^2\right)\right|_{S^{2m+1}} = \left.\left(\textstyle\sum_j dz_j\wedge d\bar z_j\right)\right|_{S^{2m+1}}.
\end{equation}
We represent $S^1=\{e^{2\pi i\theta}|\theta\in\RR\}$, so $\op{Lie}(S^1)^*=\RR \frac1{2\pi}\partial_\theta^*$ ($2\pi\partial_\theta$ is an integral coordinate).
In this basis, a moment map for the Hamiltonian diagonal action of $S^1$ on $\CC^{m+1}$ with respect to the symplectic form $\frac{i}{2\pi}\sum_j dz_j\wedge d\bar z_j$ 
is 
$$z\mapsto 1-\|z\|^2$$ 
(by setting the constant to $1$), cf. \cite[\S2.3]{dS01}. 
In particular, by \eqref{FSform},\eqref{FSisstandard},  $\omega_{\op{FS}}$ is the symplectic reduction of the form $\frac{i}{2\pi}\sum_jdz_j\wedge d\bar z_j$ on $\CC^{m+1}$.

In order to obtain the desired moment map for $\omega$, one may proceed as in \cite[\S6.6]{dS01} as follows. 
The $(S^1)^n$ action induced by $\varphi$ on $\CC^{m+1}$ has moment map
$$(w_0,...,w_{m})\mapsto -\sum_{j=0}^{m} |w_j|^2 a_j$$
with respect to $\frac{i}{2\pi}\sum_j dw_j\wedge d\bar w_j$ and a Lie algebra basis as above, cf. \cite[Exc.9]{dS01}. 
The diagonal $S^1$ action commutes with the $(S^1)^n$ action and one can take successive symplectic reductions.
One deduces that the moment map of the natural $(S^1)^n$ action on $(\CC^*)^n$ with respect to $\omega$ is
\begin{equation}
\label{momentomega}
f_\omega:(\CC^*)^n\ra\RR^n, \qquad w\mapsto -\frac{\sum_{j=0}^m |\varphi_j(w)|^2 a_j}{\sum_{j=0}^m |\varphi_j(w)|^2},
\end{equation}
see \cite[\S6.6]{dS01}, cf. \cite[\S4.2]{Fu93}.
In particular, if we are given a projective toric variety $\PP_\Delta$ containing $(\CC^*)^n$ as a dense orbit and given by a lattice polytope 
$\Delta\subset \RR^n$, we may choose the $a_j$ as the set of vertices of $\Delta$ which turns $\varphi$ into the restriction of the rational map $\PP_\Delta\ra\PP^m$ induced by linear system of $\shO_{\PP_\Delta}(1)$ with the basis of characters $\{z^{a_j}|a_j \hbox{ is a vertex of }\Delta\}$. 
We denote the resulting map by $\varphi_\Delta$ and $\omega_\Delta$ denotes the symplectic form obtained from the $\varphi_\Delta$ by pulling back $\omega_{\op{FS}}$ as above.
We have $\im f_\omega = -\Int(\Delta)$ by \cite[\S4.2]{Fu93} which is bounded unlike in Example~\ref{Cstarn}. 
Since the complex manifold underlying the mirror is $\shT_{-\Int(\Delta)}/\Lambda$, we have
\begin{center}
\fbox{
\emph{the mirror dual of $((\CC^*)^n,\Omega,\omega_\Delta)$ is a poly-annulus with cross-section $\exp(2\pi\Int(\Delta))$},}
\end{center}
\noindent see also \cite[Prop. 4.2]{Au07}.
We obtain the potential $K$ relating $\omega_\Delta$ and $\Omega$ most easily by comparing \eqref{dualviaK} and \eqref{FSform}, i.e., solving
$$ 
2i\partial\bar\partial (K\circ f_\Omega) = \frac{i}{2\pi}\partial\bar\partial\log\sum_{j=0}^m|\varphi_j|^2
$$
for $K$ which yields
$$
K( y_1,\dots, y_n)=\frac{1}{4\pi}\log\big(\sum_{j=0}^m \varphi_j(e^{-2\pi y_1},\dots,e^{-2\pi y_n})^2\big).
$$
Altenatively, we could solve the system $\check y_i = \partial_{y_i}K$ where $y_i,\check y_i$ are as in \eqref{integrals}. 
We know $y_i=(f_\omega)_i$ from \eqref{momentomega} and $\check y_i=\frac{-1}{2\pi}\log|w_i|$ from Example~\ref{Cstarn}.
Checking back the above K, we find that indeed
\begin{align*}
\partial_{y_i} K( y_1,\dots, y_n)&=
\frac{\sum_{j=0}^m \partial_{ y_i}(\varphi_j(e^{-2\pi y_1},\dots,e^{-2\pi y_n})^2)}{4\pi\sum_{j=0}^m \varphi_j(e^{-2\pi y_1},\dots,e^{-2\pi y_n})^2}\\
&=
-\frac{\sum_{j=0}^m (e^{-4\pi\sum_{k=0}^m a_{jk} y_k}) a_{ji}}{\sum_{j=0}^m \varphi_j(e^{-2\pi y_1},\dots,e^{-2\pi y_n})^2}
=
f_\omega(e^{-2\pi y_1},\dots,e^{-2\pi y_n})_i.
\end{align*}
The boundedness of $\im f_\omega$ is reflected in the asymptotic behaviour of the potential towards infinity, see e.g., Figure~\ref{figpotentials} on the right.
\end{example}

Figure~\ref{figpotentials} shows the potentials for the construction of the mirror of $\C^*$ in Example~\ref{Cstarn} and \ref{furtherCstarn} respectively. In the latter case $\omega$ is obtained via the map $\varphi:\C^*\ra(\C^*)^2, \varphi=\{1\}\times\id_{\C^*}$, i.e., $a_{01}=0,a_{11}=1$. This corresponds to taking $\Delta=[0,1]$ which is also the closure of all tangent slopes to $K$.
Let us dwell on this for a moment and motivate the next section. Consider the sequence of symplectic forms on $\CC^*$ given by $\omega_{r\Delta}$ for $r\in\NN$. This corresponds to taking the sequence of embeddings $\varphi^r_{\Delta}$ inducing a sequence of potentials $K_r( y)=\frac1{4\pi} \log(1+e^{-4\pi r y})=K(r y)$ whose normalization has the limit
\begin{equation}
\label{discretizepotential}
\lim_{r\to\infty} \frac1r K(r y)= \left\{\begin{array}{ll}
-y&\quad\text{for } y\le 0\\
 0&\quad\text{for } y\ge 0.
\end{array}\right.
\end{equation}
Thus, looking at Figure~\ref{figpotentials} on the right, the sequence of potentials approaches the piecewise linear function indicated by the positive real axis and the dotted line. Using this piecewise linear function, one can give a discrete version of the Legendre transform as we do in the following section.

\section{Large volume and large complex structure limit}
Theoretical physicists studied Calabi-Yau manifolds in order to construct conformal field theories. To obtain such a theory from the more general concept of a quantum field theory (also via a Calabi-Yau manifold), a certain function needs to vanish (the $\beta$-function, see \cite{Clay09}, \S3.2.6.2) which can be enforced by taking a \emph{large volume limit}. Since mirror symmetry is really about conformal field theories (at least by its origin), taking certain limits is an important step for its understanding. 
There are two related types of limits we are supposed to take, namely referring to \eqref{integrals},
\begin{equation} \label{limits}
\begin{array}{cclcl}
\int_{\Gamma_i}\omega&\to&\infty& &\hbox{ \emph{large volume limit,}}\\
\int_{\Gamma^*_i}\Im\Omega&\to&\infty& &\hbox{ \emph{large complex structure limit.}}
\end{array}
\end{equation}
Both of these limits amount to rescaling the affine base manifold $B$. 
Note that these interchange under mirror symmetry: a large volume limit on $X$ turns into a large complex structure limit on $\check X$ and vice versa.

We intend to take both limits simultaneously. One needs to be a bit careful about how this works with the right choice of a potential: let us first rescale the coordinate $y$ in (B) by $r$ and see how this changes everything. All data become $r$-dependant which we indicate by making $r$ an index. 
We set $y_{r,j}=ry_{j}$ and have 
$$z_{r,j}=rz_j, \quad\partial_{y_{r,j}}=\frac1r\partial_{y_j},\quad \Lambda_r=\frac1r\Lambda\quad\hbox{and}\quad\Omega_r=r^n\Omega.$$
The potential is as before determined by $\omega$ and this in turn is determined by the condition that the integral over a path scales by $r$: a priori, there are different ways to obtain an $r$-dependant potential:
\begin{enumerate}
\item The first option is to just take the pullback of $K$ via $y_{r,j}=ry_{j}$. This is $K_r'(y_r):=K(\frac1r y_{r})=K(y)$. In terms of dual coordinates, this leads to 
$\check y'_j(y_r)=\partial_{y_{r,j}}K_r'(y_r)=\frac1r\partial_{y_j}K(y)=\frac1r \check y_j(\frac1r y_{r})$. This is not what we want because it means that while enlarging the $y$-coordinates, we shrink the $\check y$-coordinates.
\item The next option is pulling back the dual coordinates via $y_{r,j}=ry_{j}$, i.e., set $\check y'_{r,j}=\check y_j(\frac1r y)$. With the previous calculation, it is easy to see that this corresponds to taking for the new potential the scaled pullback
$K^0_r(y_r):=rK'_r(y_r)$.
\item Finally, in order to actually take the large volume limit simultaneously as the large complex structure limit, we need to scale the dual coordinates as well, i.e., $\check y_{r,j}= r\check y_{j}$. This is realized by rescaling the pullback potential even more by taking $K_r(y_r):=r^2K'_r(y_r)$.
\end{enumerate}
There are two types of limits that typically occur: metric limits and algebraic limits, for a discussion, see \cite[7.3.6.]{Clay09}.
In some sense, these are represented by the two potentials shown in Fig. \ref{figpotentials}.
Note that if we choose $K(y)=y^2/2$ then  $K_r(y_r)=y_r^2/2$, so this potential remains invariant under taking the simultaneous limit. 
The effect is that the base $B$ of $f_\omega$ and $f_\Omega$ becomes longer and longer as one approaches the limit. 
Rescaling the metric to normalize the diameter yields $B$ itself as a limit the Calabi-Yaus. 
For an elliptic curve with potential $y_r^2/2$, the metric limit is thus a circle, cf. \cite{Gr08},\,Conj.\,5.4.
We are interested in algebraic limits and for such, the non-self-dual second potential in Figure~\ref{figpotentials} is more relevant.
\begin{figure} 
\resizebox{0.8\textwidth}{!}{
\input{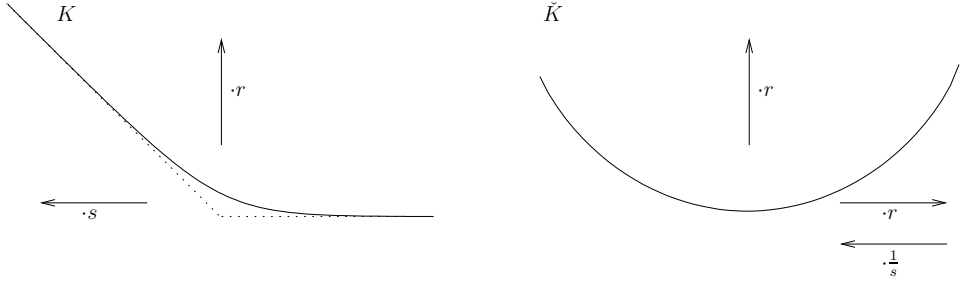} 
}
\caption{$K(y)=\frac1{4\pi} \op{log}(1 + e^{-4\pi y})$ and $\check K(y)=\frac1{4\pi}((y+1)\op{log}(y+1)-y\op{log}(-y))$}
\label{scalings}
\end{figure}
Figure~\ref{scalings} illustrates how the scaling of the potential (here by factor $r$) and the base coordinate (here by factor $s$) influences the Legendre dual and dual potential. The diagram really only shows part of all rescaling options where the remaining ones come from applying the given ones on the dual side. In fact, in view of Figure~\ref{scalings}, the result of scaling by $r$ on either side results
in scaling both potentials by $r^2$ and both coordinates by $r$ as we did in 3. above. 
There is still the degree of freedom of scaling by $s$ which has a reciprocal effect on the dual. This explains why
the limit we gave in \eqref{discretizepotential} appears to be turned into the limit $r\to 0$ now. 
In truth, it was a limit with respect to the parameter $s$. The important point is that in algebraic examples, there is a non-trivial rescaling by $s$ but it is non-homogenous along the base, i.e., in some regions it looks like a contraction, in others like an expansion. We will see this in the algebraic degeneration of an elliptic curve as well as in the mirror duality of $\PP^1$. It is really this rescaling which yields a discretization of the Legendre transform. Before we give an example, we relate $r$ to the algebraic coordinate:
on $\shT_B/\frac1r\Lambda$, we consider the two potentials derived from $K_r,K_r^0$ but modified by some inhomogeneously rescaling by some $s$.
These potentials lead to symplectic forms $\omega,\omega_0$ via the first equation in \eqref{dualviaK} and we have $\omega=r\omega_0$. 
Reparametrizing $|t|=e^{-2\pi r}$ for $t$ a coordinate on the unit disk, we get
$$\omega_t = \frac{\log |t|}{-2\pi}\,\omega_0.$$
While $\omega_t$ is going to infinity as $t\to 0$, we will find that $\omega_{0}$ is bounded. 

\begin{example}[Elliptic curve]
\label{example-elliptic-curve}
The elliptic curve has been considered from an SYZ perspective many times before. We mostly follow \cite{Gr08},\,\S6, see also \cite{Clay09}, \S8.4.1: We fix $n\in\NN$ and consider the affine manifold $B=\RR/n\ZZ$ with $y$ being the standard coordinate on $\RR$ and obtain the elliptic curve $X_r=\shT_B/ \frac1{r}\Lambda$ with periods $1$ and $irn$. 
The family parameter $r$ can be complexified: either ad hoc by using the complex coordinate $t$ on the unit disc as before and then $X_t=\shT_B/\frac{\log (t^n)}{2\pi i}\Lambda$ (note that we abuse notation here, we use the identification $\shT_B=\CC/n\ZZ$ via (B)) or more conceptually by invoking the $B$-field as in \cite{Clay09}, \S6.2.3. The limit for $t\ra 0$ can be filled by a cycle of $\PP^1$s of length $n$. This turns the total space of the family into a maximally unipotent degeneration\footnote{This means $\shX$ is flat over the base such that $X=\shX_{t_0}$ for some $t_0\neq 0$ and $T\in\End(H^\bullet(X,\QQ))$, the monodromy operator around the special fibre at $t=0$, satisfies $(T-\id)^{n+1}=0$ and $(T-\id)^{n+1}\neq 0$ with $n=\dim X$.}. 
Siebert had the idea to use log geometry to view the singular special fibre $\shX_0$ as a (log) smooth Calabi-Yau. 
Indeed, let us compute the logarithmic cotangent sheaf on $X_0$, i.e., the restriction of the relative logarithmic cotangent sheaf 
$K_{X_0}=\Omega^1_{X_0}(\log X_0):=\Omega^1_{\shX/O}(\log \shX_0)|_{\shX_0}$ with $O=\hbox{unit disk}$. For each irreducible component $\PP^1$ of $\shX_0$, we have $K_{\shX_0}|_{\PP^1}=\Omega^1_{\PP^1}(\log(\{0\}\cup\{\infty\}))\cong \shO_{\PP^1}$ and locally at an intersection point the pair $(\shX,\shX_0)$ is $(\Spec \CC[u,v],V(t))$ with $t=uv$ and thus 
$$\Omega^1_{\CC^2/\CC_t}(\log V(t))=(\shO_{\CC^2}\frac{du}{u}\oplus \shO_{\CC^2}\frac{dv}{v})/ \shO_{\CC^2}(\frac{du}{u}+\frac{dv}{v})\cong\shO_{\CC^2}.$$ We deduce $K_{\shX_0}\cong\shO_{\shX_0}$, so $\shX_0$ is a log elliptic curve. To obtain a nowhere vanishing global section $\Omega$ of $K_{\shX_0}$ we can just extend the local section $\frac{du}{u}$ in a standard chart of one of the components.
There is a (degenerate) Strominger-Yau-Zaslow fibration $\shX_0\ra B$ given as the compactification of the special Lagrangian fibration (with respect to $\Omega$ and $\omega_0$) on the dense subset of $\shX_0$ whose intersection with each $\PP^1$ is $\CC^*=\PP^1\backslash(\{0\}\cup\{\infty\})$.

Let us discuss the potential $K$ and K\"ahler form. We already mentioned that $K=y^2/2$ is not a useful choice here. In fact, Gross realized \cite[\S6]{Gr08}, that if we take an open cover of $X_0$ in $\shX$, the intersection of the nearby fibre with a neighbourhood of a node approaches $\shT_{(0,1)}/\Lambda_r$ for $r\to\infty$ whereas away from the nodes it approaches $\shT_{[0,0]}/\Lambda_r$, so all the mass in the complex geometry goes to the nodes. Conversely, all the mass in the symplectic geometry should leave any small neighbourhood of any node.
This of course depends on the choice of potential which we make as follows. 

In general, we want to have a relatively ample line bundle $\shL$ on $\shX$ and sections $s_0,...,s_m$ which are in bijection with the zero-dimensional strata (which are the nodes of $X_0$ in this example) $v_0,...,v_m$ in $X_0$ whose vanishing locus is contained in $X_0$ and such that $s_j$ vanishes along precisely those components of $X_0$ that do not contain $v_j$. 
In analogy to Example~\ref{furtherCstarn}, we then define the family of two-forms $\omega_r=\frac{i}{2\pi}\partial\bar\partial\log\sum_j|s_j|^{2r}$ on $\shX$ that is fibrewise a symplectic form.
Let $\omega$ denote the two-form on $\shX\setminus X_0$ that restricts to $\omega_{\frac{\log|t|}{-2\pi}}$ on the fibre $X_t$. Its normalization is $\omega_0 = \frac{-2\pi}{\log|t|}\omega$. 

In our example, this limit is the potential given by \eqref{discretizepotential} on each $\PP^1$ component of $X_0$ (up to the addition of an affine function). Indeed, only two $s_j$ are non-vanishing on this $\PP^1$ and they give the potential on the right of Figure \ref{figpotentials}. 
For concreteness, let us refine the example by considering the family of Fermat elliptic curves in $\PP^2$ given by $z_0z_1z_2+t(z_0^3+z_1^3+z_2^3)=0$, then $\shX$ is the blow-up of $\PP^2$ in the base locus of the family, $X_0=\{z_0z_1z_2=0\}$, $B=\RR/3\ZZ$, $\shL$ can be chosen as $\shO(1)$ and $s_j=z_j$. 
The upshot is: the potential on $B$ becomes a piecewise affine function with non-linearity at the three integral points of $B/3\ZZ$ corresponding to the equators of the components of $X_0$.
\end{example}

The example led us to the consideration of a piecewise affine potential in the limit. 
We deal with a version of the Legendre transform for such potentials in the next section.
Moreover, so far we have been dealing only with the situation where all fibres of the SYZ maps are smooth. Talking about compact manifolds with vanishing first Chern class, this restricts one to the study of complex tori, e.g., the elliptic curve just studied. For Hyperk\"ahler manifolds or Calabi-Yau manifolds in the strong sense\footnote{This means $h^\bullet(X,\shO_X)=h^\bullet(S^n,\QQ)$ for $n=\dim X$.}, one has to allow singular torus fibres. The critical loci of these fibres play an important role in the theory. There is another way of obtaining interesting geometry, namely by allowing a boundary for the affine manifold over which the SYZ fibration takes lower-dimensional tori as fibres, a typical situation for the compactifications of moment maps from $(\CC^*)^n$.

\section{Algebraic limits and the discrete Legendre transform of a tropical manifold}
We have already seen in an example that an algebraic large complex structure limit with simultaneous large volume limit leads to a discretization of the Legendre transform. 
\begin{figure}
\resizebox{0.7\textwidth}{!}{
\input{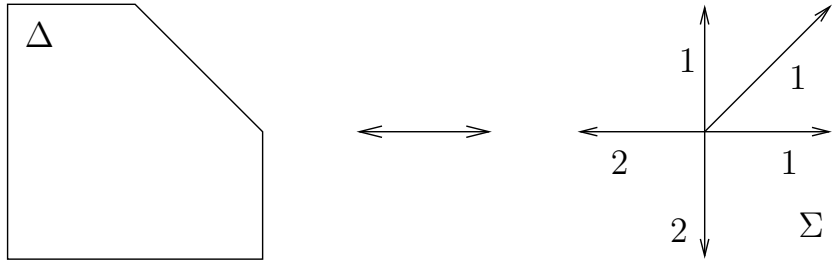}
}
\caption{A polytope is dual to a fan with piecewise linear function. The piecewise linear function is given up to addition of a linear function by its slope changes along the rays as given in the diagram}
\label{poly_dualto_fan}
\end{figure}
A general definition of this has been given in \cite{GS06} for an affine manifold with (a certain type of) singularities that behave well with regard to the piecewise affine potential. 
A discrete Legendre transform on a vector space had been known before, see \cite[\S14]{Ar78}.
We are going to give a natural extension to manifolds with polyhedral boundary. The simplest example of a discrete Legendre transform is the correspondence 
$$\Delta \leftrightarrow (\Sigma,\varphi)$$
of a polytope with a fan and piecewise linear convex function, well-known in toric geometry, see \cite{Fu93},\,\S3.4 as well as Figure~\ref{poly_dualto_fan}. 
The underlying manifolds are the polytope $\Delta$ and a real vector space respectively.
Note that $\varphi$ now plays the role of the strictly convex function $\check K$, but we need to weaken the assumption on $K,\check K$ from strictly convex as in the smooth case to just convex.\footnote{Confusingly in the discrete world (e.g. \cite{Fu93}), for a piecewise affine function on a polyhedral complex the notion \emph{strictly convex} is used for the property where the maximal cells coincide with non-extendable domains of linearity of the function. This is actually the type of function we want.}

The definition of a piecewise linear function $\varphi$ associated to a polytope $\Delta$ can be given as
\begin{equation}
\label{varphitoric}
\varphi(n) = \max\{ \langle n,m \rangle \,|\, m\in\Delta \}
\end{equation}
where $\langle\cdot,\cdot\rangle$  denotes the pairing of a vector space with its dual space.
If we take $K\equiv 0$ for the piecewise linear function on the polytope,
this coincides with the previous definition since one can check (e.g. \cite[\S14]{Ar78}) that
$$\check K(\check y) = \max_y\{\sum_i\check y_i y_i - K(y)\}.$$
Note that \cite{Fu93} uses $\varphi(n) =-\inf\{ \langle n,m \rangle \,|\, m\in\Delta \}$. 
We should make a remark on sign conventions here that also explains the minus sign in \eqref{momentomega}.
Since our discussion is governed by the Legendre transform and this associates to a point the tangent at a convex function over the point, positive directions should get mapped to positive directions under this transform unlike in \cite{Fu93} where concave functions are used.

The general construction of a discrete Legendre transform is obtained from patching this example in both directions: assume we have an integral affine manifold $B$, i.e., a real affine manifold with an atlas whose transition functions are in $\ZZ^n\rtimes\op{GL}_n(\ZZ)$. Moreover, we assume to have a polyhedral decomposition $\P$ of $B$, i.e., $\P$ is a set of lattice polytopes each of which comes with an immersion in $B$, the set $\P$ covers $B$, is closed under intersection in $B$ and two polytopes in $\P$ coincide if there image in $B$ does. We also need a polarization $\varphi$ which is a section of $\op{PAC}(B,\RR)/\op{Aff}(B,\RR)$, the sheaf of piecewise affine convex functions on $B$ (piecewise with respect to $\P$) with rational slopes modulo the sheaf of affine functions on $B$ (both with rational slopes). We require that the non-extendable domains of linearity of $\varphi$ coincide with the maximal cells in $\P$. We also require the boundary of $B$ to be locally convex, more precisely, near each point in $\partial B$, the pair $(B,\partial B)$ looks like an open subset of a lattice polytope with its boundary. Such a triple $(B,\P,\varphi)$ is called a \emph{tropical manifold}.
The \emph{discrete Legendre transform} (DLT) associates another tropical affine manifold to $(B,\P,\varphi)$ and is a duality:
$$(B,\P,\varphi) \longleftrightarrow (\check B,\check\P,\check\varphi).$$
The dual is constructed as follows: The neighbourhood of each vertex $v$ in $\P$ can be identified with a neighbourhood of the origin of a fan $\Sigma_v$ and $\varphi$ restricts to a piecewise linear convex function on its support. Thus from the duality in Figure~\ref{poly_dualto_fan}, we obtain a lattice polytope $\check v$. On the other hand, for each maximal cell $\sigma$ in $\P$, again by the duality in Figure~\ref{poly_dualto_fan}, we obtain a fan 
$\check \Sigma_\sigma$ with a piecewise linear function $\check\varphi_\sigma$. 
\begin{figure}
\resizebox{0.9\textwidth}{!}{
\input{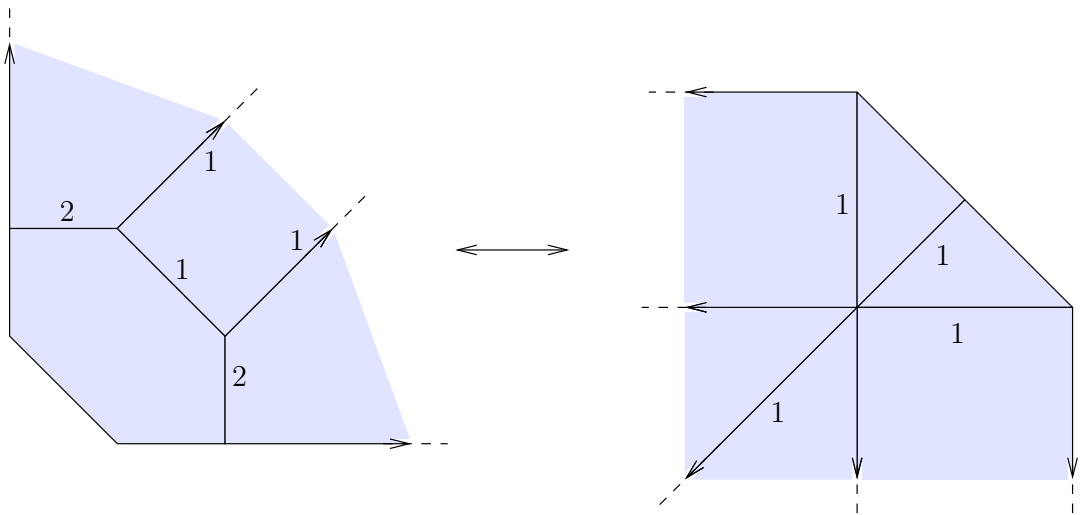}
}
\caption{An example of a discrete Legendre transform}
\label{DLT}
\end{figure}
Finally, $\check B$ is given by gluing all these polytopes and fans according to their adjacency, see Figure~\ref{DLT} for an example.

\begin{example}[Duality of cones as a DLT]
\label{dualityofcones}
Note that the duality of cones is a special case of a DLT: Let $\sigma\subset \RR^n$ be a rationally generated polyhedral cone containing no non-trivial linear subspace and
$$\check\sigma=\{n\in\Hom(\RR^n,\RR)\,|\,n(m)\ge 0 \hbox{ for all }m\in\sigma\}.$$
Taking trivial piecewise linear functions and for the polyhedral decompositions the set of faces respectively gives a discrete Legendre transform
$$\sigma\longleftrightarrow-\check\sigma.$$
Note that this is more general that the polytope-to-fan duality (e.g., Figure~\ref{poly_dualto_fan}) because given a polytope $\Delta$, we may take $\sigma$ to be
$$\Cone(-\Delta) =\{(rm,r)\,|\,m\in-\Delta,r\in\RR_{\ge0}\}\subset \RR^n\times\RR,$$
the cone over $-\Delta$. Then $(\Sigma,\varphi)$, the DLT of $\Delta$, is obtained from the dual cone $\check\sigma\subseteq \Hom(\RR^n\oplus\RR,\RR)$ as follows: $\Sigma$ is the projection of the proper faces of $\check\sigma$ under the restriction
$$\Hom(\RR^n\oplus\RR,\RR)\sra \Hom(\RR^n,\RR)$$
and the graph of $\varphi$ is the section of this projection given by $\partial\check\sigma$.
\end{example}
We will come back to this example later.

It is quite remarkable that the construction of the discrete Legendre transform even works if the tropical manifold has singularities as long as the local monodromy around the singularities respects the polyhedral decomposition, see \cite{GS06}. The discriminant loci in $B$ and $\check B$ are then homeomorphic.
Singularities are an important feature of the story. If one wants to study compact Calabi-Yau manifolds, the base $B$ of the SYZ fibration needs to be a homology sphere, see \cite[Prop. 2.37]{GS06}. Therefore, the fibration needs to have singular fibres which are reflected in the base as singular locus of the affine structure of codimension two. The singularities only affect the affine structure, the underlying topological space will still be a topological manifold (with boundary).
The local monodromy on the tangent bundle $\shT_B$ around a branch of the discriminant coincides with the monodromy in the cohomology of a nearby smooth torus fibre for case (A) and the homology for (B).
See \cite{Gr05} for a systematic account on how to obtain a DLT from reflexive polytopes and nef partitions. For further examples on affine manifolds with singularities, see \cite{HS03}, \cite{Zh98}, \cite{CM06} and \cite{Ru05}.


Before closing this section, we would like to introduce natural refinements of the cell decompositions $(B,\P)$ and $(\check B,\check\P)$ 
that give topologically a common refinement on the interiors of $B$ and $\check B$.
This is given by the \emph{barycentric subdivision} 
(cf. \cite[Def. 1.25]{GS06} for the compact case and \cite[Def. 3.2]{Ts13} for an alternative definition in the non-compact case with the draw-down that doesn't seem natural in the context of SYZ fibrations). 
The definition we give requires that each unbounded cell $\tau\in\P$ has the property that the convex hull of its vertices $\op{conv}(\tau^{[0]})$ is a face of $\tau$.
This is satisfied in Fig.~\ref{DLT}.

We define a triangulation $\P^{\op{bar}}$ of $B$, 
which introduces one new vertex in each relative interior of a compact cell $\tau\in\P$. This vertex is the \emph{barycenter} of the cell and is defined as the average of the cell's vertices 
$v^{\op{bar}}_\sigma=\frac1{\#\sigma^{[0]}}\sum_{v\in \sigma^{[0]}} v$ where $\sigma^{[0]}$ denotes the set of vertices of $\sigma$. 
We may use the same definition to associate a barycenter to an unbounded cell, so for $\tau$ unbounded we have $v^{\op{bar}}_\tau=v^{\op{bar}}_{\op{conv}\tau^{[0]}}$.
We then set
$$\P^{\op{bar}} = \{ \op{conv}\{ v^{\op{bar}}_{\tau_0},\dots,v^{\op{bar}}_{\tau_k}\}\,|\,\tau_0\subsetneq\dots\subsetneq\tau_k,\tau_i\in\P,k\ge 0\}\cup \P^{\op{bar}}_{\op{unbounded}}$$
where $\op{conv}$ means taking the convex hull and $\P^{\op{bar}}_{\op{unbounded}}$ will be empty if each cell in $B$ is bounded. It is defined as 
$$\P^{\op{bar}}_{\op{unbounded}} = \{ \op{conv}\{ v^{\op{bar}}_{\tau_0},\dots,v^{\op{bar}}_{\tau_k}\}
+\sum_{i=1}^k\RR_{\ge0}\rho_{\tau_i}
\,|\,\tau_0\subsetneq\dots\subsetneq\tau_k,\tau_i\in\P,\tau_i\hbox{ is unbounded}\}$$
where $\rho_{\tau_i}$ is the sum of all primitive integral generators of the rays in $\tau_i$, so $\sum_{i=1}^k\RR_{\ge0}\rho_{\tau_i}$ is a cone generated by such rays and its sum with 
$\op{conv}\{ v^{\op{bar}}_{\tau_0},\dots,v^{\op{bar}}_{\tau_k}\}$ should be read as a Minkowski sum (i.e., pointwise sum).

Note that indeed $\P^{\op{bar}}$ is a refinement of the polyhedral decomposition $\P$ of $B$ and, after removing the boundary respectively, topologically also of $\check \P$ of $\check B$, namely by respectively merging all cells in $\P^{\op{bar}}$ which contain a vertex that is in $\P$ but not in $\P^{\op{bar}}$. 
\begin{figure}
\begin{minipage}[r]{0.5\textwidth}
\input{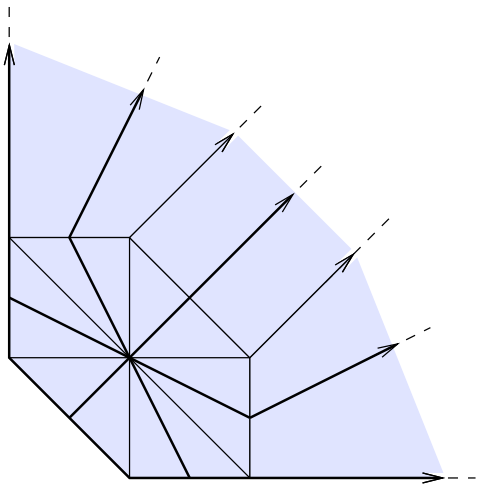}
\end{minipage}
\begin{minipage}[l]{0.4\textwidth}
\caption{The barycentric subdivision of the left hand side in Figure~\ref{DLT} indicating by bold lines how to obtain the cell decomposition given by its right hand side}
\label{bary}
\end{minipage}
\end{figure}
See Figure~\ref{bary} for an example.

\section{The degenerate Calabi-Yau fibre and the reconstruction problem}
\label{section_reconstruct}
As in the case of the elliptic curve, the tropical manifold $(B,\P,\varphi)$ encodes a degenerate fibre as follows (\cite{GS06}): 
each cell $\sigma\in\P$ gives a projective toric variety $$\PP_\sigma = \Proj\CC[\Cone(\sigma)\cap(\ZZ^n\oplus\ZZ)]$$
where 
$\Cone(\sigma)$
was defined in Example~\ref{dualityofcones}.
This is functorial for inclusions of cells: $\tau\subseteq\sigma\Rightarrow \PP_{\tau}\subseteq\PP_\sigma$, so we may form the limit
$$\check X_0(B,\P,\varphi) := \varinjlim_{\sigma\in\P}\PP_\sigma$$
which is called the degenerate Calabi-Yau in the \emph{cone picture}. This should be thought of as a degeneration of (A) in \S\ref{section1}. Dually, concerning a degeneration of (B), for each $\sigma\in\P$, we may consider the \emph{fan along $\sigma$} by which we mean the following. Let $U_\sigma$ be a sufficiently small neighbourhood of the relative interior of $\sigma$. The image of $\{\tau\in\P\,|\,\sigma\subset\tau\}$ under the projection 
$U_\sigma\sra U_\sigma /\sigma$ (where two points are identified if their difference is parallel to $\sigma$) gives a neighbourhood of the origin of a fan $\Sigma_\sigma$ in $\RR^{n-\dim\sigma}$ unique up to isomorphism. Let $X_{\Sigma_\sigma}$ denote the corresponding toric variety. This construction is contravariantly functorial for inclusions: $\tau\subseteq\sigma\Rightarrow X_{\Sigma_\sigma}\subseteq X_{\Sigma_\tau}$, so we may form
$$X_0(B,\P,\varphi) := \varprojlim_{\sigma\in\P}X_{\Sigma_\sigma}$$
which we call the degenerate Calabi-Yau in the \emph{fan picture}. It is not hard to see that in fact
$$X_0(B,\P,\varphi) =\check X_0(\check B,\check\P,\check\varphi)$$
which should be compared to Figure~\ref{logandmoment}:
indeed, there is a continuous map
$$f_\omega: \check X_0(B,\P,\varphi)\ra B$$
by taking a direct limit over all moment maps $f_{\omega_\sigma}:\PP_\sigma\ra \sigma$ for each $\sigma\in\P$, see Example~\ref{furtherCstarn} for the definition of $\omega_\sigma$. 
This does not coincide with the limit map $f_{\omega_0}$ discussed in Example~\ref{example-elliptic-curve} but it is $f_{\omega_1}$ restricted to the central fibre.
The meaning of $f_\omega$ could be understood as this: suppose we have a nearby fibre $X_t$, then we can use symplectic parallel transport to get a retraction map $\check X_t\ra \check X_0$ and we can compose this with $f_\omega$ to get a Lagrangian fibration $\check X_t\ra B$. It is currently not clear how to turn this into a special Lagrangian fibration. We have a diagram:

$$
\xymatrix@C=30pt
{
\check X_0(B,\P,\varphi) \ar[d]^{\check f_{\omega}} &
X_0(\check B,\check\P,\check\varphi) \ar[d]^{f_{\omega}}\\
(B,\P,\varphi) \ar@{<->}[r]^{\op{DLT}}& (\check B,\check\P,\check\varphi)
}
$$
We have called $X_0(B,\P,\varphi)$ and $\check X_0(B,\P,\varphi)$ \emph{Calabi-Yau}. This is justified if its canonical bundle is trivial. These spaces have a log structure would be entirely encoded in 
$\P$ for the first and in $\varphi$ for the second if were no singularities in the affine structure. 
The singularities however contribute non-discrete moduli of the log structure encoded in so-called \emph{slab functions}, see \cite{GS11}. We will not go into defining log structures, but recall that in the case of the elliptic curve we constructed a sheaf of log differential forms which was trivial. This generalizes as long as the transition functions of $B$ can be chosen in $\ZZ^n\rtimes\op{SL}_n(\ZZ)$, i.e., $B$ is orientable.
The log differential forms restricted to each component 
$\PP_\sigma$ of $\check X_0(B,\P,\varphi)$ are just
$\Omega^k_{\PP_\sigma}(\log D_\sigma),$
the differential forms with logarithmic poles along $D_\sigma$
where $D_\sigma$ is the complement of the dense torus in $\PP_\sigma$. 
These sheaves glue to a sheaf
$\Omega^k:=\Omega^k_{\check X_0(B,\P,\varphi)^\dag/\Spec\CC^\dag}$,
though the gluing is non-trivial whenever singularities appear (the dagger indicating the presence of a log structure), see \cite[\S3.2]{GS10}.
If $B$ is orientable, $\Omega^n$ is trivial and a section gives a global holomorphic volume form with logarithmic poles. 

The \emph{reconstruction problem} is the question of whether one can reconstruct a smooth (or at most orbifold) Calabi-Yau $X_t$ from its degeneration $X_0$; more precisely, whether we can lift $X_0$ from a space over a point to a flat family $\shX$ over the unit disk whose  non-zero fibres have at most orbifold singularities.
In general, so in presence of singularities, this is a very difficult problem towards which Gross and Siebert accomplished a major break-through in \cite{GS11} by proving a canonical liftability to 
$\Spec \CC\llbracket t\rrbracket$ assuming that the local monodromy of the affine singularities of $B$ cannot be factored (locally rigid). The parametrization of the disk is also important and Gross and Siebert obtain the one trivializing the Gauss-Manin connection (flat coordinates). Their proof is constructive and involves \emph{wall-crossings}. We will come back to this in a later section.

We now treat an easy case: Assume that $B$ is a lattice polyhedron in $\RR^n$ and $\P$ a subdivision of it given by a piecewise linear function $\varphi$. So in particular, we have no singularities. It is not hard to see that the dual $(\check B,\check\P,\check\varphi)$ also has the property that it globally embeds in a vector space (the dual space). The DLT here can be worked out as follows.
Let $\Delta(B,\P,\varphi)$ be the polyhedron in $\RR^n\oplus\RR$ given as
$$\Delta_{(B,\P,\varphi)} = \{(m,r)\in \RR^n\oplus\RR| \varphi(m)\ge r\}$$
and let 
$$\Sigma_{(B,\P,\varphi)} =\{0\}\cup \{\overline{\Cone(\tau)}\,|\,\tau\in\P\}
$$
be the fan in $\RR^n\oplus\RR$ 
where $\overline{\Cone(\tau)}$ denotes the closure of $\Cone(\tau)$ in $\RR^n\oplus\RR$. We define the piecewise linear function
$\varphi_{(B,\P,\varphi)}(m,r)=r\varphi(m)$. We have a DLT
$$\Delta_{(B,\P,\varphi)}\leftrightarrow 
(\Sigma_{(\check B,\check\P,\check\varphi)},\varphi_{(\check B,\check\P,\check\varphi)})$$
which is really just the classical toric story as in Figure~\ref{poly_dualto_fan}.
The original DLT $(B,\P,\varphi)\leftrightarrow (\check B,\check\P,\check\varphi)$ is contained in this as a ``sub-DLT'' by intersecting with $\RR^n\times\{1\}$.
Moreover, this picture solves the reconstruction problem:
The fan $\Sigma_{( B,\P,\varphi)}$ maps to the fan of $\AA^1$ by the projection to the second factor $\RR^n\oplus\RR\sra\RR$, so we have a map of toric varieties
$$f:\shX( B,\P,\varphi):=X_{\Sigma_{( B,\P,\varphi)}}\ra \Spec\CC[t]$$
such that $f^{-1}(0)=X_0( B,\P,\varphi)$ and
$f^{-1}(t)$ is irreducible for $t\neq 0$. In fact
$f^{-1}(t)$ is isomorphic to the toric variety given by the \emph{asymptotic fan} of $(B,\P,\varphi)$ which is just the sub-fan of $\Sigma_{( B,\P,\varphi)}$ contained in $\RR^n\times\{0\}$. So if this gives a smooth toric variety, a general fibre of $f$ is smooth.
This is the total space description for the fan picture. 

There is a dual version, the cone picture $\check\shX( \check B,\check\P,\check\varphi)$ of the total space satisfying 
$$\check\shX( \check B,\check \P,\check \varphi)=\shX( B,\P,\varphi).$$ 
Gross and Siebert use this cone picture description to prove the more general reconstruction (non-embedded situation). Those familiar with toric geometry will know that we have
$$\check\shX(B,\P,\varphi)=\Proj \CC[\overline{\Cone(\Delta_{( B,\P,\varphi)})}\cap\ZZ^{n+2}].$$
Let us recall how this works by gluing charts: To each vertex $v$ of $\Delta_{(B,\P,\varphi)}$ we associate the ring $R_v=\CC[\RR_{\ge0}(\Delta_{(B,\P,\varphi)}-v)\cap(\ZZ^n\oplus\ZZ)]$ which is naturally a $\CC[t]$-algebra by mapping $t$ to the monomial given by the unique generator of the second summand in $\ZZ^n\oplus\ZZ$ (indeed, it is contained in $\Delta_{(B,\P,\varphi)}-v$). The affine varieties $\Spec R_v$ will give an open cover of $\shX$. The intersection of two such,
$\Spec R_v$ and $\Spec R_w$, is empty if no cell in $\P$ contains both $v$ and $w$ and otherwise for $\tau$ being the minimal cell containing both, we may localize the rings $R_v$ and $R_w$ by inverting all elements that are sums of monomials with exponents contained in $\RR_{\ge0}(\tau-v)$ (respectively $\RR_{\ge0}(\tau-w)$).
Denoting the resulting rings $R_{v,\tau}$ and $R_{w,\tau}$, we have a natural isomorphism $R_{v,\tau}\ra R_{w,\tau}$ induced by 
$$\RR_{\ge0}(\Delta_{(B,\P,\varphi)}-v)+\RR(\tau-v)
=\RR_{\ge0}(\Delta_{(B,\P,\varphi)}-w)+\RR(\tau-w).$$
All these isomorphisms are compatible and glue to give $\check\shX(B,\P,\varphi)$ and a map $\check f:\check\shX(B,\P,\varphi)\ra\CC$ such that $\check f^{-1}(0)=\check X_0(B,\P,\varphi)$. To see the latter, note that we identify 
\begin{equation}
\label{lifttoDelta}
\Proj \CC[\Cone(\sigma)\cap(\ZZ^n\oplus\ZZ)] \cong
\Proj \CC[\Cone(\varphi(\sigma))\cap((\ZZ^n\oplus\ZZ)\oplus\ZZ)].
\end{equation}

\section{Compactifying divisors and the Landau-Ginzburg potential}
We have already dealt with the situation where $B$ has a boundary when we discussed discrete Legendre transforms. We now want to match it with the discussion of SYZ fibrations from \S\ref{section1}. For this, let us consider the mirror dual of $\PP^1$. We have already treated the mirror dual of $\CC^*$ with respect to its Fubini-Study-metric coming from the embedding in $\PP^1$, we have
$$\RR\stackrel{f_\Omega}{\longleftarrow}\CC^*
\stackrel{f_\omega}{\longrightarrow} (0,1).$$
The map $f_\omega$ naturally extends to $\PP^1\ra[0,1]$. We may think of the compactifying divisor $D=\{0\}\cup\{\infty\}$
as adding (partially) contracted SYZ fibres. In fact, we contract the 1-cycle which we used to define our base coordinate via the first integral in \eqref{integrals}. Phrased differently, the holomorphic cylinders\footnote{It is possible to choose them holomorphic, in fact there is a natural choice.} which we used to define the base coordinate $y$ on $(0,1)$ becomes a holomorphic disk. By the maximum principle, there are no holomorphic disks in $\CC^*$, but they do appear as we compactify to $\PP^1$. 
It is insightful to interpret the presence of holomorphic disks from the point of view of Floer theory, see \cite{Au09} for a detailed account. We already mentioned that the mirror $\check X$ of $X=\CC^*$
can be considered as the moduli space of pairs $(L,\nabla)$ where $L$ special Lagrangian tori with a $U(1)$-connection $\nabla$ on $L\times\CC$. Fukaya-Oh-Ohta-Ono \cite{FO$^3$08} give an obstruction for the intersection Floer homology complex to be a complex. If we are interested in the Floer homology $HF^\bullet(\shL,\shL)$ of $\shL=(L,\nabla)$ with itself, the obstruction is
\begin{equation}
\label{m0}
m_0(\shL)=\sum_{{\beta\in\pi_2(X,L)}\atop{\mu(\beta)=2}} n_\beta(\shL) z_\beta(\shL)
\end{equation}
where $n_\beta(\shL)$ is the (virtual) number of holomorphic disks of homotopy class $\beta$ which contain a pre-determined general marked point in $L$, $\mu(\beta)$ denotes the Maslov index of $\beta$ and 
\begin{equation}
\label{potentialmonomialsasintegrals}
z_\beta(\shL)=\exp(-\int_\beta\omega)\op{hol}_\nabla(\partial\beta)\in\CC^*
\end{equation}
for $\op{hol}_\nabla(\partial\beta)$ the holonomy of $\nabla$ along 
$\partial\beta$.
The important observation is that $z_\beta(\shL)$ gives a holomorphic function on $\check X$. Just note its similarity with the holomorphic coordinate
$$w_j=\exp(2\pi i(x_j+i\int_{\Gamma_j}\omega))$$ on $\check X$ given in \S\ref{section1}.
By \cite{Au07}, Lemma 3.1, the condition $\mu(\beta)=2$ is equivalent to $\beta.D=1$ where $D$ is the compactifying divisor and the dot denotes the algebraic intersection number. We learn that a partial compactification of $X$ yields a holomorphic function $m_0$ on $\check X$ (assuming that $\eqref{m0}$ has finitely many summands or converges).
Motivated by physics, this function is called a \emph{Landau-Ginzburg-potential} (LG potential) and denoted $W$. 
The pair $(\check X, W)$ is called a \emph{Landau-Ginzburg model} (LG model). 
In fact more generally, an LG model will simply be a variety with a flat holomorphic function to $\CC$ as well as a restriction of such to an open subset in the analytic topology.
\begin{figure}
\input{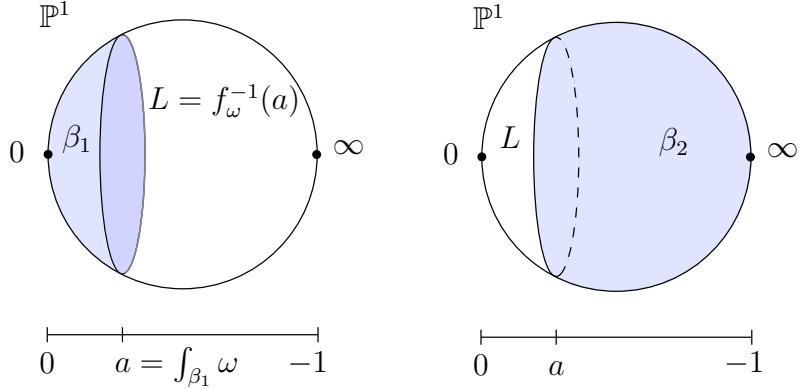}
\caption{The two holomorphic disks giving the LG potential of the mirror of $\PP^1$}
\label{P1}
\end{figure}
Coming back to the example of $\PP^1$, the two-point compactification of $\check X$, we obtain a LG potential on 
$$X=\{w\in\CC^*\,|\, e^{2\pi 0}<|w|<e^{2\pi 1}\}$$
given by $W=e^{2\pi}(w+\frac1w)$, see Figure~\ref{P1}. We could have gotten rid of the factor $e^{2\pi}$ had we rescaled $m_0$. 
This generalizes to smooth toric Fano varieties, see \cite{Au09}~Prop.2.5:
\begin{proposition} 
\label{Fanoprop}
Given a smooth projective toric Fano variety $\PP_\Delta$, the LG potential on its mirror is given by
$$W=\sum_{\tau\subset\Delta\hbox{\,{\tiny is a facet}}} e^{-2\pi\alpha_\tau}z^{n_\tau}$$
where $n_\tau\in\Hom(\RR^n,\RR)$ is the primitive integer inward normal vector to $\tau$, such that $\tau$ is given by intersecting the affine hyperplane $n_\tau+\alpha_\tau=0$ with $\Delta$. 
Moreover, $z^{n_\tau}$ is the character associated to $n_\tau$ for the torus containing the poly-annulus which is the mirror of the dense $(\CC^*)^n$ in $\PP_\Delta$.
\end{proposition}
Note that we may also study non-compact Fanos, e.g., by embedding
$\CC^n$ in $\PP^n$, we have that the mirror dual of $\CC^n$ is the LG model
$$\{(w_1,...,w_n)\in(\CC^*)^n\,|\,{1}<|w_i|<e^{2\pi}\}\stackrel{w_1+...+w_n}{\lra}\CC.$$
Let us now consider the large volume limit of this picture.
By taking $\lim_{r\to\infty}r\omega$, we enlarge the mirror poly-annulus until it becomes all of $(\CC^*)^n$. The potential will also move to infinity, but can be normalized similarly as we normalized the symplectic form previously, see \cite{Au07},\,\S4.2. Under normalization, it remains the same and we have in the large volume limit
\begin{center}
\fbox{
\emph{the mirror dual of $((\CC)^n,\Omega,\omega_{\PP^n})$ is the LG model $(\CC^*)^n\xrightarrow{w_1+...+w_n}\CC$}.}
\end{center}
Let us now see how we find the LG potential in the context of the discrete Legendre transform, i.e., in the degeneration limit. 
\begin{example}[LG potential on the mirror of $\PP^1$]
We consider the example of $\PP^1$ again, which is given by the cone picture
$$\PP^1=\check X_0([-1,0],\{\{-1\},\{0\},[-1,0]\},0)$$ 
and the DLT of $B=[-1,0]$ is the fan of $\PP^1$ with piecewise linear function $\varphi$ whose slope changes by $1$ at the origin. We have for the mirror degenerate Calabi-Yau the fan picture 
$$X_0:=X_0([-1,0],\{\{-1\},\{0\},[-1,0]\},0) = \AA^1\sqcup_{\{0\}} \AA^1= V(uv)\subseteq \AA^2$$
We take the potential $W_0=u+v$, i.e., the standard coordinate on each $\AA^1$.
The reconstruction of $X_0$ is given by $\shX=\AA^2\ra\AA^1, (u,v)\mapsto uv$. Let us view the same from the perspective of the cone picture. We denote the DLT of $B$ by 
$(\check B,\check\P,\check\varphi)$ and may assume $\check\varphi(0)=0$. We have the monoid algebra
$$\shX = \Spec\CC[P],\qquad P=\Delta_{(\check B,\check\P,\check\varphi)}\cap(\ZZ\oplus\ZZ).$$
Let $e_1, e_2$ be generators for the two summands of $\ZZ\oplus\ZZ$ respectively. We set $w=z^{e_1}$ and $t=z^{e_2}$. 
Since $e_2\in P$, $\CC[P]$ is a $\CC[t]$-algebra giving the map
$\Spec\CC[P]\ra \Spec\CC[t]=\AA^1.$
The generators of the $P$ are 
$e_1+\varphi(e_1)e_2$ and $-e_1+\varphi(-e_1)e_2$, so the generators of $\Spec\CC[P]$ are $wt^{\varphi(e_1)}$ and $w^{-1}t^{\varphi(-e_1)}$.
Denoting these by $u,v$, we have $\CC[P]=\CC[u,v]$.
We claim that the sum $u+v$ is the reconstruction of the LG potential:
$$W = wt^{\varphi(e_1)}+w^{-1}t^{\varphi(-e_1)}.$$
Indeed this restricts to $W_0$ on $X_0$.
Inserting the $\varphi$ as given by \eqref{varphitoric} from $\Delta=[-1,0]$, we get for $t\neq 0$
$$W = w+w^{-1}t$$
for the potential on $X_t=V(xy-t)\cong\CC^*$. Taking $t=1$ reproduces the mirror of $\PP^1$ constructed before Prop.~\ref{Fanoprop} up to a factor of $e^{2\pi}$ and up to the restriction to an annulus.
\end{example}

What we did for $\PP^1$ here generalizes directly to the case of a general $(B,\P,\varphi)$, see \cite{CPS11}. The potential $W_0$ on a component $\PP_\sigma$ of $\check X_0(B,\P,\varphi)$ 
is $0$ if $\sigma$ is compact. Otherwise, let
$\op{rays}(\sigma)$ denote the set of equivalence classes of (unbounded) extremal rays of $\sigma$ up to translation. The potential on $\PP_\sigma$ is given by
\begin{equation}
\label{rayspotential}
W_0|_{\PP_\sigma}=\sum_{(n_0+\RR_{\ge0}n)\, \in\, \op{rays}(\sigma)} z^{n}
\end{equation}
where $(n_0+\RR_{\ge0}n)$ denotes a representative of an element in $\op{rays}(\sigma)$ for which we require $n$ to be a primitive integral vector. 
Clearly $z^n$ doesn't depend on the choice of representative.
These local potentials glue to a LG potential $W_0$ on $\check X_0(B,\P,\varphi)$.
See \cite{CPS11} for a solution of the reconstruction problem for this potential. We again restrict ourselves to the easy case:
Let us assume that $(B,\P,\varphi)$ is embedded in $\RR^n$. Recall from the end of \S\ref{section_reconstruct} the local description of the total space 
$\check\shX=\check\shX(B,\P,\varphi)$ of the smoothing of $\check X_0(B,\P,\varphi)$.
For each vertex $v\in\P$, we have an affine chart $\Spec R_v$ of $\check\shX$. The reconstructed potential 
$W:\check\shX\ra\CC$ is given in each $R_v$ by the sum
\begin{equation}
\label{DLTpotential}
W = \sum_{
{(n_0+\RR_{\ge0}n)\, \in\, \bigcup\{\op{rays}(\sigma)|\sigma\in\P\}}}
z^n t^{\varphi(n+n_0)-\varphi(n_0)}.
\end{equation}
It can be shown that $\varphi(n+n_0)-\varphi(n_0)$ is an invariant of the equivalence class of $n_0+\RR_{\ge0}n$.
Note that this indeed restricts to $W_0$ on $X_0$
making use of the identification \eqref{lifttoDelta}.
It is also in line with the above example for the mirror of $\PP^1$ where the sum was $u+v$ and we had $\varphi(n_0)=\varphi(0)=0$.
More generally in the presence of singularities of the affine structure, one needs to sum over all \emph{broken lines} which we have implicitly done here, too. Broken lines are an analogue of holomorphic disks in tropical geometry. See \cite{Gr09},\cite{CPS11} for more details.

\section{Mirror duality for Landau-Ginzburg models}
\label{MirdualLGmodels}
We are now in the position to study a duality of Landau Ginzburg models.
We understood in the first section that the mirror dual of $(\CC^*)^n$ is again
$(\CC^*)^n$ or some analytic open subset thereof depending on the choice of symplectic form. We understood in the previous section that partial compactifications on one side lead to a LG potential on the other side. LG models are well-known to be the mirror duals of projective Fano varieties, some of which are compactifications of $(\CC^*)^n$, some others (possibly all) can be degenerated torically such that the mirror is also obtained from the given discrete Legendre transform construction. However, in principle, there is nothing stopping us from looking at partial compactifications of $(\CC^*)^n$ on both sides as in Fig.~\ref{DLT}, e.g., the reader will meanwhile hopefully agree with the slogan
\begin{center}
\fbox{
\emph{the mirror dual of\quad $\CC^n\xrightarrow{w_1+...+w_n}\CC$\quad is
\quad $\CC^n\xrightarrow{w_1+...+w_n}\CC$.}}
\end{center}
The discrete Legendre transform underlying this slogan is the duality of very simple cones, namely
$$\RR^n_{\ge 0} \stackrel{\op{DLT}}{\longleftrightarrow} -\RR^n_{\ge 0},$$
more precisely, one $\RR^n_{\ge 0}$ sits in the dual space of the vector space containing the other. While the DLT provides a very general framework for the construction of very sophisticated Landau-Ginzburg models (e.g., with singularities in the affine structure), we give here a simple and yet very useful subset of the wide range of DLT duals:

Let us fix a free abelian group $M\cong \ZZ^n$, $M_{\RR}=M\otimes_{\ZZ}\RR$, $N=\Hom_{\ZZ}(M,\ZZ)$, $N_{\RR}=N\otimes_{\ZZ}\RR$.
Consider a strictly\footnote{This means it doesn't contain a non-trivial linear subspace.} convex rational polyhedral cone $\sigma\subseteq M_{\RR}$
with $\dim\sigma=\dim M_{\RR}$, and let $\check\sigma\subseteq N_{\RR}$ be
the dual cone, 
\[
\check\sigma:=\{n \in N_{\RR}\,|\,\hbox{$\langle n,m\rangle \ge 0$ for all
$m\in\sigma$}\}. 
\]
We already explained in Example~\ref{dualityofcones} that the duality
$\sigma\leftrightarrow-\check\sigma$ constitutes a DLT. For the simplicity of the exposition, we remove the minus sign from $\check\sigma$ in the following and call $\sigma\leftrightarrow\check\sigma$ and related constructions a DLT.
Note that in this notation, the previous slogan results from starting with the cone $\sigma= \RR_{\ge 0} e_1\oplus...\oplus \RR_{\ge 0}e_n$ where $e_1,...,e_n$ is a basis of $M$. Note that if $e_1,...,e_n$ were only a basis of $M\otimes_\ZZ\QQ$ but not of $M$, we would already be studying an interesting duality of quotient singularities (in fact this relates to the Berglund-H\"ubsch construction \cite{BH92}), cf. \cite{Bo10}. Let us remain in the smooth world. So since the corresponding toric varieties
\begin{align*}
\check X_0(\check\sigma)=X_0(\sigma)=X_{\sigma} & =\Spec \CC[\check\sigma\cap N]\\
\check X_0(\sigma)=X_0(\check\sigma)=X_{\check\sigma} & =\Spec \CC[\sigma\cap M]
\end{align*}
are usually singular, we choose toric desingularizations by choosing
fans $\Sigma$ and $\check\Sigma$ which are refinements
of $\sigma$ and $\check\sigma$
respectively, with $\Sigma$ and $\check\Sigma$ consisting only of standard
cones, i.e., cones generated by part of a basis for $M$ or $N$.

We now obtain smooth toric varieties $X_{\Sigma}$ and $X_{\check\Sigma}$.
However, the resolution has broken the DLT property: $\Sigma$ is not the DLT of $\check\Sigma$ in general. This can be fixed as follows. We may assume that we have chosen resolutions given by a piecewise linear functions $\varphi$, $\check\varphi$ respectively. Then there are polytopes 
$P\subseteq M_\RR$, $\check P\subseteq N_\RR$ such that we have DLTs
\begin{equation}
\label{middleDLTs}
\begin{array} {rcl}
(\sigma,\Sigma,\varphi)&\leftrightarrow& \check P\\
(\check\sigma,\check \Sigma,\check \varphi)&\leftrightarrow& P.
\end{array}
\end{equation}
Moreover, these have the property that
$$\begin{array} {rcl}
\overline{\Cone(P)}\cap M_\RR&=&\sigma\\
\overline{\Cone(\check P)}\cap N_\RR&=&\check\sigma,
\end{array}$$
where the overline means taking the closure and the cones are contained in
$M_\RR\oplus\RR$ (resp. $N_\RR\oplus\RR$) so that intersection with $M_\RR$
(resp. $N_\RR$) makes sense.
Note that we have the fan pictures
$X_\Sigma=X_0(\sigma,\Sigma,\varphi)$, 
$X_{\check\Sigma}=X_0(\check\sigma,\check\Sigma,\check\varphi)$.
By the construction in the previous section, we obtain reconstructed potentials
$\check W:\check\shX(\sigma,\Sigma,\varphi)\ra\CC$, $ W:\check\shX(\check\sigma,\check \Sigma,\check \varphi)\ra\CC$, which make sense to write down as elements 
$$
\begin{array}{lcl}
\check W=\sum_{{\RR_{\ge0}n \hbox{\tiny\ is a ray in }\check\Sigma}\atop{n\in N\hbox{\tiny\  is primitive}}}
z^n t^{\check\varphi(n)}
&\in& \CC[\overline{\Cone(\check\sigma)}\cap (N\oplus\ZZ)]
=\CC[\check\sigma\cap N]\otimes_\CC\CC[t]\\
W=\sum_{{\RR_{\ge0}m \hbox{\tiny\ is a ray in }\Sigma}\atop{m\in M\hbox{\tiny\  is primitive}}}
z^m t^{\varphi(m)}
&\in&\CC[\overline{\Cone(\sigma)}\cap (M\oplus\ZZ)]
=\CC[\sigma\cap M]\otimes_\CC\CC[t].
\end{array}
$$
So the potentials pull back from $X_\sigma\times\AA^1_t,X_{\check\sigma}\times\AA^1_t$ respectively. For a fixed $t$, we have diagrams
\begin{equation} \label{mirrorLGs}
\xymatrix@C=30pt
{
 & X_{\Sigma}\ar[d]\ar[ld]_W&& X_{\check\Sigma}\ar[d]\ar[rd]^{\check W}\\
\CC & X_{\sigma}\ar[l]^W&& X_{\check\sigma}\ar[r]_{\check W}&\CC\\
}
\end{equation}
One needs to take a close look to observe that this duality is actually ``balanced'' in the following sense. One might wonder what happens if one chooses a different resolution $\Sigma_{\op{new}}$ instead of $\Sigma$. Then $X_\Sigma$ becomes $X_{\Sigma_{\op{\new}}}$ but the potential $W$ remains ``the same'' (being the pullback of the same potential on $X_\sigma$). However, while on the dual side $X_{\check\Sigma}$ remains the same space, its potential $\check W$ changes to $\check W_{\op{new}}$ because it is a sum over all rays in $\Sigma_{\op{\new}}$. So it is not possible to change only one side by adding exceptional divisors.
Of course, the geometry of $X_{\Sigma_{\op{\new}}}$ might be very different from that of $X_{\Sigma}$, e.g., one of them might have a trivial canonical bundle while the other has a more positive one. To ensure that the geometry of 
$X_\Sigma$ doesn't differ considerably from that of $X_\sigma$, we would want
$X_\Sigma\ra X_\sigma$ to be a crepant resolution. Such does not always exist in the category of smooth schemes, however it does exist in general in the category of orbifolds which should be the slightly more general framework to be used here.

\begin{figure}
\begin{center}
\resizebox{0.82\textwidth}{!}{
\ \ \input{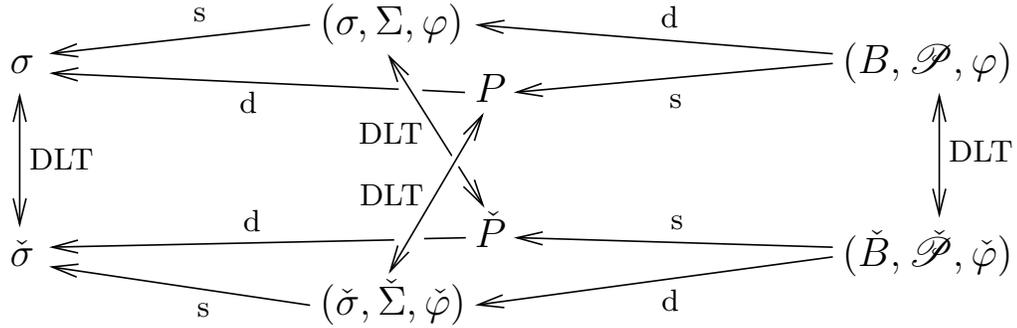}
}
\end{center}
\caption{Tropical manifolds and their relationships: {\scriptsize DLT}
marks a discrete Legendre transform (up to sign convention),  {\scriptsize s} marks a
subdivision,  {\scriptsize d} marks a deformation/degeneration}
\label{DLTs}
\end{figure}
While the balancing argument just given is a weak one to rectify mirror symmetry, we should actually argue by the discrete Legendre transform.
There are four DLTs in place three of which we have seen already, 
see Figure~\ref{DLTs}. It has been shown in \cite{GKR12} that there exists a (non-unique) DLT pair $(B,\P,\varphi)\leftrightarrow (\check B,\check\P,\check\varphi)$
which ``dominates'' the two DLTs given in \eqref{middleDLTs}.
Most importantly, the potentials constructed for 
$\shX(B,\P,\varphi)$, $\shX(\check B,\check\P,\check\varphi)$
via the previous section agree with $W, \check W$ respectively in the following sense: the space
$\shX(B,\P,\varphi)$ relates to $\shX(\sigma,\P,\varphi)$
by a deformation, i.e., there is a flat family with general fibre isomorphic to 
$\shX(\sigma,\P,\varphi)$ and special fibre given by $\shX(B,\P,\varphi)$. 
Moreover this family is birational to the trivial family with fibre
$\shX(\sigma,\P,\varphi)$ and the potential on $\shX(B,\P,\varphi)$ is the pullback of the potential $W$ from the trivial family.

The mirror duality of Landau-Ginzburg models given in \eqref{mirrorLGs} has been used in \cite{GKR12} to construct mirror duals for varieties which are not necessarily Fano or Calabi-Yau, e.g., for varieties of general type. A notion of mirror symmetry for such varieties didn't exist before the cited work had been started, so this relatively simple construction for duals is already quite powerful. Note also that the famous mirror construction of Batyrev-Borisov is reproducible from this duality, see \cite{GKR12}. Note that the potentials in loc.cit. had been permitted to have more general coefficients, i.e., 
$$W=\sum_{{\RR_{\ge0}m \hbox{\tiny\ is a ray in }\Sigma}\atop{m\in M\hbox{\tiny\  is primitive}}}
c_m z^m t^{\varphi(m)}$$
for some (general) $c_m\in\CC$ and similarly for $\check W$ (independently of the coefficients of $W$).
This can be argued to make sense by changing the (complexified) symplectic form on either side, recall from \eqref{potentialmonomialsasintegrals} that the monomials are integrals of the symplectic form.

\begin{figure}
\begin{center}
\resizebox{0.8\textwidth}{!}{
\input{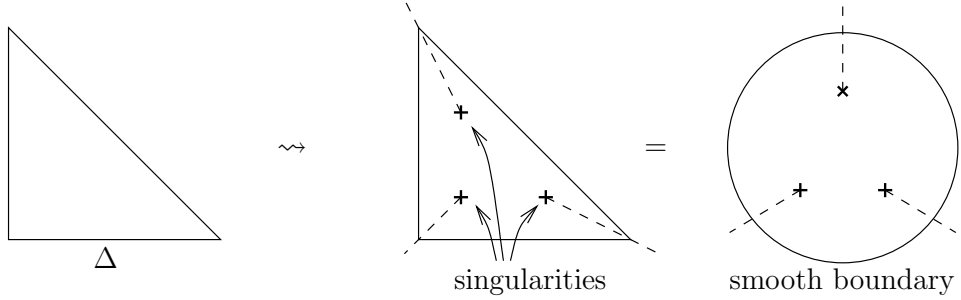}
}
\end{center}
\caption{The cone picture of $\PP^2$ and how to trade corners for singularities}
\label{pullincorners}
\end{figure}

There is yet one flaw in the picture: The potential which we give in \eqref{DLTpotential} is the ``naive potential''. It agrees with the Floer theoretic one in the Fano case by Prop.~\ref{Fanoprop}, however $X_\Sigma$, $X_{\check\Sigma}$ are rarely Fano. More generally, there will be non-rigid rational curves in $X_\Sigma$ or $X_{\check\Sigma}$ and these cause \emph{disk bubbling} and non-geometric virtual counts of holomorphic disks (see \cite{Au09}). Such give rise to (possible infinitely many) additional Maslov index two holomorphic disks and thus terms in the potential. To keep this under control, the authors of \cite{CPS11} required the boundary of $B$ and $\check B$ to be smooth. In fact they suggested to smooth the boundary by trading ``corners'' in $B$ (or $\check B$) for singularities of the affine structure of $B$ (or $\check B$), see Figure~\ref{pullincorners}.

The advantage is that the tropical potential (the generalization of \eqref{DLTpotential} to affine manifolds with singularities) for a smooth boundary of $B$ (or $\check B$) seems to agree with the Floer theoretic one. The additional terms arise from holomorphic disks attaching to the singularities in the SYZ fibration and these can be accounted for tropically. We shall study this for an example in the next section. 
Let us record here the main result of \cite{GKR12} which supports the mirror duality \eqref{mirrorLGs} from a cohomological point of view. For this, the general setup of \eqref{mirrorLGs} is restricted to the situation where
$\sigma$ has the special shape of a Gorenstein cone, i.e., there is a lattice polytope $\Delta$ such that 
$$\sigma=\Cone(\Delta).$$
For this to make sense, we need to write $M_\RR$ as $(M_0\oplus\ZZ)\otimes_\ZZ\RR$ where $M_0\cong \ZZ^{n-1}$ and 
$\Delta\subset M_0\otimes_\ZZ\RR$. Note that the existence of a toric crepant resolution $X_\Sigma\ra X_\sigma$ is equivalent with the existence of a triangulation $\P$ of $\Delta$ into simplices for which the edges emanating from a vertex in each form a basis of $M_0$.
The authors of \cite{GKR12} prove the following:
\begin{theorem}
Assume that $\Delta$ has at least one interior lattice point, $\PP_\Delta=\Proj\CC[\sigma\cap M]$ is smooth and that there is a projective crepant toric resolution $X_\Sigma\ra X_\sigma$ factoring through the blowup of the origin $X_\Sigma\ra\Bl_0 X_\sigma\ra X_\sigma$ then the blow-up of the origin $X_{\check\Sigma}=\Bl_0 X_{\check\sigma}\ra X_{\check\sigma}$ is a toric resolution. The diagram \eqref{mirrorLGs} specializes to
$$
\xymatrix@C=30pt
{
 & X_{\Sigma}\ar[d]^{\hbox{\tiny crepant}}\ar[ld]_W&& \Tot(\shO_{\PP_\Delta}(-1))\ar[d]\ar[rd]^{\check W}\\
\CC & X_{\sigma}\ar[l]^W&& X_{\check\sigma}\ar[r]_{\check W}&\CC.\\
}
$$
where $\Tot(\shL)=\Spec (\Sym(\shL^{-1}))$ denotes the total space of a line bundle.
The critical locus of $\check W$ is a hypersurface $S\subset \PP_\Delta$ which is smooth if the coefficients of $\check W$ were chosen general. The Kodaira dimension of $S$ is
$$\kappa(S) = \min\{\dim \Delta',n-2\}$$
where $\Delta'$ is the convex hull of the lattice points in the interior of $\Delta$. We have that 
$$W^{-1}(0)=D_{v_1}\cup ... \cup D_{v_r}\cup \tilde W_0$$ 
is normal crossings, $\tilde W_0$ is the strict transform of the zero fibre of $W:X_\sigma\ra\CC$ and $D_{v_i}$ are toric exceptional divisors of $X_\Sigma\ra X_\sigma$ projecting to the origin. They are indexed by the lattice points in the interior of $\Delta$. The critical set near the origin $\check S=\Sing W^{-1}(0)$ supports the sheaf of vanishing cycles $\shF_{\check S}=(\phi_{W,0}\CC)[1]$ which carries the structure of a cohomological mixed Hodge complex. Denoting
$$h^{p,q}(\check S,\shF_{\check S})=\dim\Gr^F_p\HH^{p+q}(\check S,\shF_{\check S}),$$
we have
$$h^{p,q}(S)= h^{d-p,q}(\check S,\shF_{\check S})$$
where $d=\dim S=n-2$.
\end{theorem}

\section{Moving the compactifying divisor and corrected potentials}
\begin{figure}
\resizebox{0.94\textwidth}{!}{
\input{intervalDLTs.pstex_t}
}
\caption{The four DLTs for $\Tot(\shO_{\PP^1}(-k))$}
\label{intervalDLTs}
\end{figure}
\begin{figure}
\resizebox{0.94\textwidth}{!}{
\input{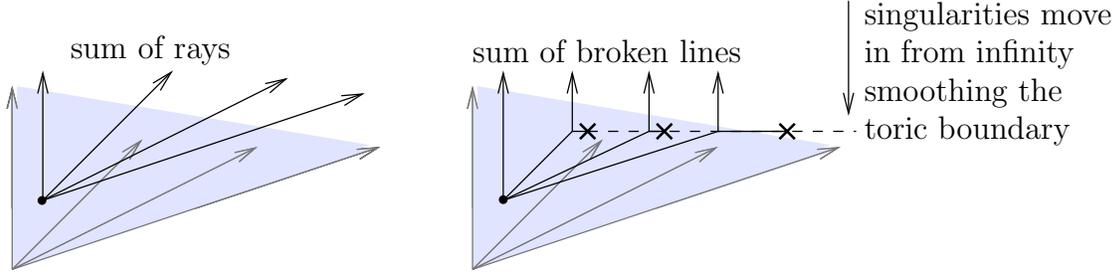}
}
\caption{Fan pictures for the minimal crepant resolution of the singularity $uv-z^3=0$ with and without smoothing of the toric boundary divisor. 
Interpreted dually, these are cone pictures for a degeneration of the singularity $C^2/\zeta_3$ where $\zeta_3$ is a primitive third root of unity acting diagonally.
The monomials in the LG potential on this singularity remain the same when smoothing the toric boundary divisor of the mirror dual: summing over rays becomes summing over broken lines}
\label{raystobrokenlines}
\end{figure}

We already mentioned the concept of trading corners for singularities, see Figure~\ref{pullincorners}. Geometrically this means the following: Recall that we started our discussion with the mirror duality of $(\CC^*)^n$ and continued by partially compactifying it to a toric variety $X_\Sigma$ using a toric divisor 
$D=X_\Sigma\backslash(\CC^*)^n$. The special Lagrangian fibration (SYZ fibration) is still entirely given on $(\CC^*)^n$ with parts of the torus fibres contracting towards $D$. There are moduli of the pair $(X_{\check \Sigma}, D)$ by moving $D$ in its equivalence class, in particular $D$ becomes non-toric by doing so. It is not known whether $X_{\check \Sigma}\backslash D$ for such a non-toric $D$ still supports a special Lagrangian fibration (using for $\Omega$ a section of 
$\Omega^n_{X_{\check \Sigma}}(\log D)$). This is already unknown for the complement of a smooth cubic in $\PP^2$. Nonetheless, we already have a good expectation of what the affine base of such a special Lagrangian fibration should look like. In the case of $\PP^2$, we depicted it on the right of Figure~\ref{pullincorners}. See \cite{Pa11} for a treatment of the case of a partial smoothing of the hyperplanes in $\PP^2$, see also \cite{CLL10}.
What happens to the mirror as we smooth $D$?
We have a natural bijection between the components of $D$ and the terms in the potential of the mirror $\check W:X_{\check \Sigma}\ra\CC$, so by smoothing $D$, we expect only one monomial to contribute to the mirror potential near $D$. On the other hand, the special Lagrangian fibration on $X_{\check \Sigma}\backslash D$ - should such exist - or at least the affine model for its base acquires singularities there are additional disks attaching to these singularities and to $D$. 
It can be checked in simple Fano examples that the monomials in the potential remain the same (up to changing coefficents) when smoothing the toric boundary divisor. 
Summing over rays in \eqref{rayspotential} is replaced by summing over \emph{broken lines} in the presence of singularities \cite{CPS11}, see Figure~\ref{raystobrokenlines}.

The singularities emanate walls (indicated dashed in Fig.~\ref{raystobrokenlines}) into the affine manifold which ought to contain the image of Maslov index zero holomorphic disks under the SYZ map $f_\Omega$ should such exist. 
These can be attached to the holomorphic disk touching $D$ and give rise to further terms in the potential. As long as $D$ itself does not contribute such walls, the tropical potential obtained in this way by counting broken lines is expected to be the correct potential meaning that it agrees with the one given in \eqref{m0}. 
Moreover the smoothing of $D$ makes $W$ proper as has been argued in \cite{CPS11}. 
The process of pulling in the corners is very ad hoc and hasn't been systematized yet. This will be treated in \cite{RS13}. 
In non-Fano cases, where Prop.~\ref{Fanoprop} possibly fails, the right count of holomorphic disks seems more accessible when the boundary divisor has been smoothed by means of counting broken lines. 
We close this article by studying the corner-pull-in-process in an example: 

\begin{example}[Corrected potential for $\Tot(\shO_{\PP^1}(-k))$ and its mirror]
\begin{figure}
\resizebox{0.94\textwidth}{!}{
\input{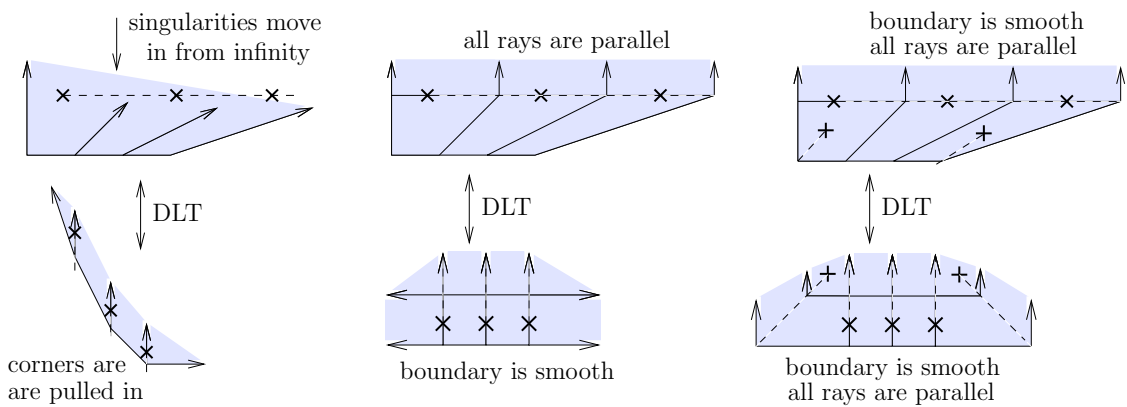}
}
\caption{Flattening the boundary}
\label{straighten}
\end{figure}

Let $\sigma=\Cone(\Delta)$ with $\Delta$ an interval of length $k$ and $\Sigma$ be the unique subdivision giving a crepant resolution of $X_\sigma$. Let $\check\sigma$ be the dual cone of $\sigma$ and $\check\Sigma$ be the fan of 
$\Tot(\shO_{\PP^1}(-k))$ which resolves $X_{\check\sigma}$. See Figure~\ref{intervalDLTs} for a how the diagram in Figure~\ref{DLTs} visualizes for this setup. We start from the DLT pair 
$(B,\P,\varphi)\leftrightarrow (\check B,\check\P,\check\varphi)$
and straighten out the boundary in these each at a time. 
See this process in Figure~\ref{straighten}. Even though we started with very simple cones, we eventually obtain a fairly interesting DLT pair whose singularities will feature scattering. The upshot of this example is that the corrections that come to the potentials don't impact the critical locus of the potential. This can be deduced from the positions of the invariant directions of the singularities towards the direction of the boundary divisor in the respective cone pictures. The critical loci together with the sheaf of vanishing cycles were the main objects of study in \cite{GKR12}.
\end{example}

\begin{acknowledgement}
The author is indebted to Bernd Siebert, Denis Auroux and Mark Gross for what he learned from them.
\end{acknowledgement}

\end{document}